\documentclass[aop,MSNbibl,nameyear,secthm,seceqn,noMRlinks,dvips]{arximspdf}
\usepackage{graphicx}

\doi{10.1214/08-AOP401}
\volume{37}
\issue{1}
\pubyear{2009}
\firstpage{206}
\lastpage{249}

\makeatletter
\newproclaim{remark}{Remark}[section]
\newtheorem{corollary}{Corollary}[section]
\makeatother

\begin{document}
\begin{frontmatter}

\title{Fractional diffusion equations and processes with randomly varying
time\thanksref{TIT1}}
\thankstext{TIT1}{Supported by ``Sapienza'' University of Rome Grant Ateneo 2007, n. 8.1.1.1.32.}
\runtitle{Processes with randomly varying
time}

\begin{aug}
\author[A]{\fnms{Enzo} \snm{Orsingher}\ead[label=e1]{enzo.orsingher@uniroma1.it}}
and
\author[A]{\fnms{Luisa} \snm{Beghin}\ead[label=e2]{luisa.beghin@uniroma1.it}}
\runauthor{E. Orsingher and L. Beghin}
\affiliation{``Sapienza'' Universit\`a di Roma}
\address[A]{Dip. di Statistica, Probabilit\`a Stat. Appl.\\
``SAPIENZA'' Universita' di Roma\\
P. Le  A. Moro 5\\
00185 Roma\\
Italy\\
\printead{e1}\\
\phantom{E-mail: }\printead*{e2}} 

\end{aug}

\received{\smonth{6} \syear{2007}}
\revised{\smonth{1} \syear{2008}}

\begin{abstract}
In this paper the solutions $u_{\nu }=u_{\nu }(x,t)$ to fractional diffusion
equations of order $0<\nu \leq 2$ are analyzed and interpreted as densities
of the composition of various types of stochastic processes.

For the fractional equations of order $\nu =\frac{1}{2^{n}}$, $n\geq 1,$ we
show that the solutions $u_{{1/2^{n}}}$ correspond to the distribution
of the $n$-times iterated Brownian motion. For these processes the
distributions of the maximum and of the sojourn time are explicitly given.
The case of fractional equations of order $\nu =\frac{2}{3^{n}}$, $n\geq 1,$ is
also investigated and related to Brownian motion and processes with
densities expressed in terms of Airy functions.

In the general case we show that $u_{\nu }$ coincides with the distribution
of Brownian motion with random time or of different processes with a
Brownian time. The interplay between the solutions $u_{\nu }$ and stable
distributions is also explored. Interesting cases involving the bilateral
exponential distribution are obtained in the limit.
\end{abstract}



\pdfkeywords{60E05, 60G52, 60J65, 33E12, 33C10, Iterated Brownian
motion, fractional derivatives, Airy functions, McKean law, Gauss--Laplace random variable, stable distributions}

\end{frontmatter}

\section{Introduction}\label{s1}

Time-fractional equations of the form
\begin{equation}\label{equation1}
\frac{\partial ^{\nu }u}{\partial t^{\nu }}=\lambda ^{2}\frac{\partial ^{2}u%
}{\partial x^{2}}, \qquad x\in \mathbb{R}, t>0,
\end{equation}%
for $0<\nu \leq 2,$ have been studied by a number of authors since the
1980s: see, for example, Wyss (\citeyear{Wys1986}), Nigmatullin (\citeyear{Nig1986}), Schneider and
Wyss (\citeyear{SchWys1989}), Mainardi (\citeyear{Mai1995}, \citeyear{Mai1996}) and, more recently, Nigmatullin
(\citeyear{Nig2006}), Angulo et al. (\citeyear{Angetal2000}, \citeyear{Angetal2005}). Hyperbolic fractional equations
similar to (\ref{equation1}) have been analyzed, for example, by Engler
(\citeyear{Eng1997}).

For exhaustive reviews on this topic, also consult Samko, Kilbas and Marichev (\citeyear{SamKilMar1993}) and
Podlubny (\citeyear{Pod1999}).

For interesting applications of fractional equations to physical problems
see, for example, Saichev and Zaslavsky (\citeyear{SaiZas1997}), Nigmatullin et al. (\citeyear{Nig2007}), Angulo
et al. (\citeyear{Angetal2005}).

Fractional diffusion equations of order $1\leq \nu <2$ emerge in the study
of the distribution of the local time of pseudoprocesses related to
higher-order heat-type equations; see Beghin and Orsingher (\citeyear{BegOrs2005}).

The time-fractional derivative appearing in (\ref{equation1}) must be
understood in the sense of Dzerbayshan--Caputo, that is%
\begin{eqnarray*}
&&\frac{\partial ^{\nu }u}{\partial t^{\nu }}(x,t)\\
&&\qquad=\cases{
\displaystyle\frac{1}{\Gamma (m-\nu )}\int_{0}^{t}\frac{1}{(t-s)^{1+\nu -m}}\frac{%
\partial ^{m}u}{\partial t^{m}}(x,s)\,ds,&\quad for $m-1<\nu <m$, \cr
\displaystyle\frac{\partial ^{m}u}{\partial t^{m}}(x,t),&\quad for $\nu=m$,
}
\end{eqnarray*}
where $m-1=\lfloor \nu \rfloor $.

Considering the derivative in the sense of Dzerbayshan--Caputo permits us to
study initial value problems for (\ref{equation1}) with initial data
represented by derivatives of integer order; on this topic, consult Mainardi (\citeyear{Mai1996}).

We assume, in particular, the following initial condition:
\begin{equation}
u(x,0)=\delta (x) \qquad \mbox{for }0<\nu \leq 1,  \label{condition 1}
\end{equation}%
and%
\begin{equation} \label{cond1}
\cases{
u(x,0)=\delta (x), \cr
u_{t}(x,0)=0,
} \qquad \mbox{for }1<\nu \leq 2.
\end{equation}

The general solution to equation (\ref{equation1}) subject to (\ref%
{condition 1}) or (\ref{cond1}) is well known [see Podlubny (\citeyear{Pod1999}), formula
(4.22), page 142] and reads
\begin{eqnarray} \label{equation2}
u_{\nu }(x,t) &=&\frac{1}{2\lambda t^{\nu /2}}\sum_{k=0}^{\infty }\frac{%
( -|x|/(\lambda t^{\nu /2})) ^{k}}{k!\Gamma (-\nu k/2+1-\nu
/2)}\nonumber\\[-8pt]\\[-8pt]
&=&\frac{1}{2\lambda t^{\nu /2}}W_{-\nu /2,1-\nu /2}\biggl( -%
\frac{|x|}{\lambda t^{\nu /2}}\biggr)  \nonumber
\end{eqnarray}%
where $W_{\alpha ,\beta }$ in (\ref{equation2}) denotes the so-called Wright
function, whose general form is
\begin{equation} \label{Wright}
W_{\alpha ,\beta }(x)=\sum_{k=0}^{\infty }\frac{x^{k}}{k!\Gamma (\alpha
k+\beta )},\qquad \alpha >-1,\beta >0,x\in \mathbb{R}.
\end{equation}

Some properties of the Wright function are investigated in Mainardi and
Tomirotti (\citeyear{MaiTom1995}) and in Gorenflo, Mainardi and Srivastava (\citeyear{GorMaiSri1998}). Initial value problems (as
well as problems on half-lines with boundary conditions) for equations like (%
\ref{equation1}) are extensively treated and solved in Mainardi
(\citeyear{Mai1994}, \citeyear{Mai1995}, \citeyear{Mai1995b}), Gorenflo and Mainardi (\citeyear{GorMai1997}) and Buckwar and Luchko
(\citeyear{BucLuc1998}).

It has been proved also that $u_{\nu }$ is nonnegative and integrates to
one for all $0<\nu \leq 2$; see, for example, Orsingher and Beghin (\citeyear{OrsBeg2004}).

We present here some alternative forms of the solution $u_{\nu }$ of (\ref%
{equation1}), either as integral functions like%
\[
u_{\nu }(x,t)=\frac{1}{\pi \nu |x|}\int_{0}^{+\infty }e^{-w}e^{-
|x|w^{\nu /2}/(\lambda t^{\nu /2})\cos ( \nu \pi /2)
}\sin \biggl( \frac{|x|w^{\nu /2}}{\lambda t^{\nu /2}}\sin \biggl(\frac{\nu \pi }{2}%
\biggr)\biggr) \,dw,
\]
or in terms of stable densities
\[
p_{\alpha }(x;\gamma ,\eta )=\frac{1}{2\pi }\int_{-\infty }^{+\infty
}e^{-i\beta x}\exp \bigl\{ -\eta |\beta |^{\alpha }e^{-i\pi \gamma /2
\beta /|\beta |}\bigr\} \,d\beta ,\qquad \alpha \neq 1,
\]
as%
\[
u_{\nu }(x,t)=\cases{
\displaystyle\frac{1}{\nu |x|^{2/\nu +1}}p_{\nu /2}\biggl( \frac{1}{|x|^{%
2/\nu }};\frac{\nu }{2},\frac{1}{\lambda t^{\nu /2}}\biggr),
&\quad $0<\nu \leq 1$, \cr
\displaystyle\frac{1}{\nu }p_{2/\nu }\biggl( |x|;\frac{2}{\nu }(\nu -1),\lambda ^{%
2/\nu }t\biggr),& \quad $1\leq \nu <2$.
}
\]

In Orsingher and Beghin (\citeyear{OrsBeg2004}), we proved that in the special case $%
\nu =\frac{1}{2},$ the solution (\ref{equation2}) coincides with the
distribution of the process%
\begin{equation}\label{iterated}
\mathcal{I}_{1}(t)=B_{1}(|B_{2}(t)|),\qquad t>0,
\end{equation}%
called the iterated Brownian motion, which consists of a Brownian motion $%
B_{1}$ whose ``time'' is an independent
reflecting Brownian motion.

In Beghin and Orsingher (\citeyear{BegOrs2003}) we have generalized this result to the case
where $\nu =\frac{1}{n},$ $n\in \mathbb{N}$. In this case, for $\lambda
^{2}=1/2$, the solution (\ref{equation2}) coincides with the distribution of
the process
\begin{equation}\label{iterated Brownian motion}
J_{1/n}(t)=B_{1}\Biggl( \prod_{j=1}^{n-1}G_{j}(t)\Biggr)
,\qquad n>1, t>0,
\end{equation}%
where the vector process $(G_{1}(t),\ldots,G_{n-1}(t))$ has the following joint
distribution:
\begin{eqnarray}\label{mult}
p(w_{1},\ldots,w_{n-1})=\frac{n^{(n-1)/2}}{(2\pi )^{(n-1)/2}\sqrt{t}%
}e^{-(w_{1}^{n}+\cdots+w_{n-1}^{n})/\sqrt[n-1]{n^{n}t}}w_{2}\cdots
w_{n-1}^{n-2},\nonumber\\[-8pt]\\[-8pt]
\eqntext{w_{j}\geq 0, 1\leq j\leq n-1,}
\end{eqnarray}%
for $n\geq 2$.

In (\ref{iterated Brownian motion}) the role of ``time'' is played by the product of independent,
positive-valued r.v.s, which cannot be identified with well-known
distributions as in the special case (\ref{iterated}).

In the special case $n=2$, we note that $\emph{J}_{1/2}(t)=\mathcal{I%
}_{1}(t)$, because (\ref{mult}) becomes the distribution of a reflecting
Brownian motion.

We are now able to prove a much stronger result for the case $\nu
=\frac{1}{2^{n}}$, $n\in \mathbb{N}$, and for $\lambda
^{2}=2^{1/2^{n}-2}$, which has a number of interesting
consequences. We will show below that (\ref{equation2}) for $\nu
=\frac{1}{2^{n}}$ can be written down as
\begin{equation}
\label{multiple integral}
\qquad u_{1/2^{n}}(x,t)=2^{n}\int_{0}^{\infty }\cdots
\int_{0}^{\infty }\frac{e^{-x^{2}/(2z_{1})}}{\sqrt{2\pi z_{1}}}\frac{%
e^{-z_{1}^{2}/(2z_{2})}}{\sqrt{2\pi z_{2}}}\cdots \frac{e^{-z_{n}^{2}/(2t)}}{\sqrt{2\pi t}}\,dz_{1}\cdots \,dz_{n}
\end{equation}%
and this coincides with the distribution of
\begin{equation}\label{iperiterated}
\mathcal{I}_{n}(t)=B_{1}(|B_{2}(|B_{3}(\cdots
(|B_{n+1}(t)|)\cdots  )|)|),\qquad t>0,
\end{equation}%
where the $B_{j}$'s are independent Brownian motions.

The iterated Brownian motion \emph{$\mathcal{I}$}$_{1}(t)=B_{1}(|B_{2}(t)|)$
has been actively investigated and many of its properties have been obtained
by Khoshnevisan and Lewis (\citeyear{KhoLew1996}), Burdzy and San Mart\`{\i}n (\citeyear{BurMar1995}), Allouba (\citeyear{All2002}).

The connection between fractional generators of order $1/2$ and the
iterated Brownian motion $\mathcal{I}_{1}(t)$ has been studied in
Allouba and Zheng (\citeyear{AllZhe2001}) and Baeumer, Meerschaert and
Nane (\citeyear{BaeMeeNan2007}). This connection was obtained in
Orsingher and Beghin (\citeyear{OrsBeg2004}) as a particular case of
the analysis of the fractional telegraph equation.

The identity
\begin{eqnarray}\label{multiple integral bis}
&&2^{n}\int_{0}^{\infty }\cdots \int_{0}^{\infty }\frac{e^{-
x^{2}/(2z_{1})}}{\sqrt{2\pi z_{1}}}\frac{e^{-z_{1}^{2}/(2z_{2})}}{%
\sqrt{2\pi z_{2}}}\cdots \frac{e^{-z_{n}^{2}/(2t)}}{\sqrt{%
2\pi t}}\,dz_{1}\cdots dz_{n}  \nonumber\\[-8pt]\\[-8pt]
&&\qquad=\frac{1}{(2t)^{1/2^{n+1}}}\sum_{k=0}^{\infty }\frac{( -
2|x|/(2t)^{1/2^{n+1}}) ^{k}}{k!\Gamma (-k/2^{n+1}+1-1/2^{n+1})}  \nonumber
\end{eqnarray}%
shows that there is a deep connection between Wright functions and Gaussian
distributions.

For the $n$-times iterated Brownian motion $\mathcal{I}_{n}(t)$, $t>0$, we
obtain the distributions of the maximum and the sojourn time (together with
the expression of moments) and we work out in detail an explicit form of
them for the case of the classical iterated Brownian motion $\mathcal{I}%
_{1}(t),$ $t>0.$

We note that $\mathcal{I}_{n}(t)$ converges in distribution, for $%
n\rightarrow +\infty $, to a Gauss--Laplace (or bilateral exponential) random
variable, independent from $t>0$.

In Orsingher and Beghin (\citeyear{OrsBeg2004}) we have seen that for the fractional
telegraph-type equation%
\begin{equation}\label{tel}
\cases{
\displaystyle\frac{\partial u}{\partial t}+2\lambda \frac{\partial ^{1/2}u}{%
\partial t^{1/2}}=c^{2}\frac{\partial ^{2}u}{\partial x^{2}},
\cr
\displaystyle u(x,0)=\delta (x),
}\qquad x\in \mathbb{R}, t>0,
\end{equation}%
the general solution coincides with the distribution of the telegraph
process $T$ whose time is an independent reflecting Brownian motion
\begin{equation}
W(t)=T(|B(t)|),\qquad t>0.  \label{Iterated telegraph}
\end{equation}

We remark that process (\ref{Iterated telegraph}) converges to (\ref%
{iterated Brownian motion}) in the Kac sense (i.e., for $\lambda $, $%
c\rightarrow \infty $, in such a way that $\frac{c^{2}}{\lambda }\rightarrow
1)$. Related interpretations of the solutions to
\begin{equation}\label{equation3}
\cases{
\displaystyle\frac{\partial u^{2\nu }}{\partial t^{2\nu }}+2\lambda \frac{\partial ^{\nu
}u}{\partial t^{\nu }}=c^{2}\frac{\partial ^{2}u}{\partial x^{2}} ,\cr
\displaystyle u(x,0)=\delta (x)%
}
\end{equation}%
are discussed in Beghin and Orsingher (\citeyear{BegOrs2003}) and Orsingher and Beghin
(\citeyear{OrsBeg2004}). Generalized forms of the fractional telegraph equation (\ref%
{equation3}) and of its solutions can be found in Saxena, Mathai and Haubold (\citeyear{SaxMatHan2006}).

We obtain here various types of relationships between the solutions $u_{\nu
} $ for different values of $\nu .$ The first one we present is the following:%
\begin{equation}\label{intro.1}
u_{\nu }(x,t)=\frac{1}{\sqrt{\pi t}}\int_{0}^{\infty }e^{-z^{2}/(4t)}u_{2\nu }(x,z)\,dz
\end{equation}%
(valid for any $0<\nu <1)$, where $u_{2\nu }$ is the solution of (\ref%
{equation1}) with order $2\nu $ instead of $\nu .$ Formula (\ref{intro.1})
leads, for $\nu =\frac{1}{2^{n}},$ to the $n$-times iterated Brownian motion
defined in (\ref{iperiterated}), since it permits us to obtain, in an
alternative way, the relationship (\ref{multiple integral}).

In the general case, (\ref{intro.1}) shows that the process related to the
equation (\ref{equation1}) of order $\nu $ can be interpreted as the
composition of a process governed by the same equation, but with order $2\nu
,$ with a Gaussian-distributed time. We also derive the analogous
relationship%
\begin{equation} \label{intro.2}
u_{\nu }(x,t)=\int_{0}^{\infty }\frac{1}{\sqrt{4\pi \lambda w}}e^{-
x^{2}/(2\lambda w)}\overline{u}_{2\nu }(w,t)\,dw,
\end{equation}%
where%
\begin{equation}\label{intro.3}
\overline{u}_{2\nu }(w,t)=\cases{
2u_{2\nu }(w,t),&\quad $w>0$, \cr
0,&\quad $w<0$.
}
\end{equation}

Here the roles of space and time are interchanged with respect to (\ref%
{intro.1}). Therefore from (\ref{intro.2}) a further interpretation of the
solution emerges, because it coincides with the density of the process%
\[
B(\mathcal{T}_{\nu }(t)),\qquad t>0,
\]
where $B$ is a Brownian motion and $\mathcal{T}_{\nu }(t)$ is a process
independent from $B$ with a distribution for each $t$ given in (\ref%
{intro.3}).

A relationship similar to (\ref{intro.1}) and connecting $u_{\nu }$ with $%
u_{m\nu }$ is established (by applying the multiplication formula of Gamma
function) for $m\geq 3$ and $0<\nu \leq 2/m$.

Substantially different situations are encountered for the special cases $%
\nu =\frac{1}{3},$ $\nu =\frac{2}{3}$ and $\nu =\frac{4}{3}.$ In particular
for $\nu =\frac{2}{3}$ we show that the solution to (\ref{equation1})
possesses the following simple form:%
\begin{equation} \label{intro.4}
u_{2/3}(x,t)=\frac{3}{2}\frac{1}{\lambda \sqrt[3]{3t}}Ai\biggl( \frac{%
|x|}{\lambda \sqrt[3]{3t}}\biggr) ,
\end{equation}%
where $Ai(x)$ is the Airy function. The latter emerges as a solution to
third-order heat-type equations of the form%
\[
\frac{\partial u}{\partial t}=-\frac{\partial ^{3}u}{\partial x^{3}},\qquad
t>0, x\in \mathbb{R}.
\]

By using again the relationship (\ref{intro.1}) we get, for the case $\nu =%
\frac{1}{3}$, the following result:%
\begin{equation}\label{intro.5}
u_{1/3}(x,t)=\frac{3}{2}\int_{0}^{\infty }\frac{e^{-z^{2}/(4t)}%
}{\sqrt{\pi t}}\frac{1}{\lambda \sqrt[3]{3z}}Ai\biggl( \frac{|x|}{\lambda
\sqrt[3]{3z}}\biggr) \,dz.
\end{equation}

This suggests that we should interpret $u_{1/3}$ as the distribution of
\[
J_{1/3}(t)=A( |B(t)|) ,\qquad t>0,
\]
where $A$ is a process whose one-dimensional distribution is given in (\ref%
{intro.4}), which coincides with the symmetric stable process of order $1/3$.

Similar relationships seem not to hold for the solutions to fractional
equations of order $\nu =\frac{1}{n},n>3$, because the fundamental solutions to
\begin{equation}\label{n-order equation}
\frac{\partial u}{\partial t}=c_{n}\frac{\partial ^{n}u}{\partial x^{n}},
\end{equation}
$c_{n}=\pm 1$, are sign-varying functions on the whole $x$-axis (while, for $%
n=3,$ only on the negative half-line), as shown in detail in Lachal (\citeyear{Lac2003}).
Therefore they cannot be used to construct the functions $u_{\nu }$ emerging
from (\ref{equation1}), which, for $0<\nu \leq 2$, are nonnegative and
integrate to one. We note that the solutions to (\ref{n-order equation})
themselves have been represented as distributions of compositions of
artificial processes, which do not display a probabilistic structure [see
Funaki (\citeyear{Fun1979}), Hochberg and Orsingher (\citeyear{HocOrs1996}), Benanchour, Roynette and Vallois (\citeyear{BenRoyVal1999})].

Finally the previous results permit us to establish connections between the
solutions $u_{2/3^{n}}$ and $u_{2/3^{n-1}}$. Moreover the
explicit form (\ref{intro.4}) of $u_{2/3}$ suggests that we should interpret them
as distributions of processes similar to the $n$-times iterated Brownian
motion, but with the role of $B$ replaced by $A$ and the time represented by
nested products of the random variables $G_{j}$ defined in (\ref{iterated
Brownian motion}).

\section{Iterated Brownian motions generated by fractional
equations}\label{s2}

In this section we examine in detail various relationships between solutions
to diffusion equations like (\ref{equation1}) and processes involving
Brownian motion. All results of this section refer to equations of order $%
0<\nu \leq 1.$

We start with the following general theorem:

\begin{thm}\label{thm1}
The solution to
\begin{equation}
\cases{
\displaystyle\frac{\partial ^{\nu }u}{\partial t^{\nu }}=\lambda ^{2}\frac{\partial ^{2}u%
}{\partial x^{2}}, \cr
u(x,0)=\delta (x),}
\qquad x\in \mathbb{R}, t>0,  \label{eq4}
\end{equation}%
for $0<\nu \leq 1$, can be represented as
\begin{equation}\label{vologda}
u_{\nu }(x,t)=\frac{1}{\sqrt{\pi t}}\int_{0}^{\infty }e^{-z^{2}/(4t)}u_{2\nu }(x,z)\,dz
\end{equation}%
where $u_{2\nu }$ is the solution to
\begin{equation}\label{init1}
\cases{
\displaystyle\frac{\partial ^{2\nu }u}{\partial z^{2\nu }}=\lambda ^{2}\frac{\partial
^{2}u}{\partial x^{2}}, \cr
u(x,0)=\delta (x),
}\qquad \mbox{for }0<\nu \leq \frac{1}{2}
\end{equation}%
or
\begin{equation} \label{init2}
\cases{
\displaystyle\frac{\partial ^{2\nu }u}{\partial z^{2\nu }}=\lambda ^{2}\frac{\partial
^{2}u}{\partial x^{2}}, \cr
u(x,0)=\delta (x), \cr
u_{t}(x,0)=0,}
\qquad \mbox{for }\frac{1}{2}<\nu \leq 1.
\end{equation}
\end{thm}

\begin{pf}
By applying the duplication formula of the Gamma
function we have that
\begin{equation}\label{Gamma}
\Gamma \biggl( -\frac{\nu k}{2}+1-\frac{\nu }{2}\biggr) =\sqrt{\pi }2^{\nu
(k+1)}\frac{\Gamma (1-\nu (k+1))}{\Gamma ( 1/2(1-\nu
(k+1))) }.
\end{equation}%
By plugging (\ref{Gamma}) into (\ref{equation2}) we get that
\begin{eqnarray*}
u_{\nu }(x,t)
&=&\frac{1}{2\lambda t^{\nu /2}}\sum_{k=0}^{\infty }\frac{( -%
|x|/(\lambda t^{\nu /2})) ^{k}\Gamma ( 1/2(1-\nu
(k+1))) }{k!\sqrt{\pi }2^{\nu (k+1)}\Gamma (1-\nu (k+1))} \\
&=&\frac{1}{\sqrt{\pi }2^{\nu +1}\lambda t^{\nu /2}}%
\sum_{k=0}^{\infty }\frac{( -|x|/(\lambda t^{\nu /2}))
^{k}\int_{0}^{\infty }e^{-w}w^{-\nu /2(k+1)-1/2}\,dw}{k!2^{\nu
k}\Gamma (1-\nu (k+1))} \\
&=&\frac{1}{\sqrt{\pi }2^{\nu +1}\lambda t^{\nu /2}}%
\int_{0}^{\infty }e^{-w}w^{-\nu /2-1/2}\\
&&\phantom{\frac{1}{\sqrt{\pi }2^{\nu +1}\lambda t^{\nu /2}}%
\int_{0}^{\infty }}
{}\times\sum_{k=0}^{\infty }%
\frac{1}{k!\Gamma (1-\nu (k+1))}\biggl( -\frac{|x|}{\lambda 2^{\nu }(wt)^{\nu
/2}}\biggr) ^{k}\,dw \\
&=&[\mbox{in view of (\ref{equation2}) with suitable arrangements}%
] \\
&=&\frac{1}{\sqrt{\pi }2^{\nu }t^{\nu /2}}\int_{0}^{\infty
}e^{-w}w^{-\nu /2-1/2}\bigl(2\sqrt{tw}\bigr)^{\nu }u_{2\nu }\bigl(x,2\sqrt{%
tw}\bigr)\,dw \\
&=&\frac{1}{\sqrt{\pi }}\int_{0}^{\infty }e^{-w}w^{-1/2}u_{2\nu }\bigl(x,2%
\sqrt{tw}\bigr)\,dw \\
&=&\bigl[ 2\sqrt{tw}=z\bigr] \\
&=&\frac{1}{\sqrt{\pi t}}\int_{0}^{\infty }e^{-z^{2}/(4t)}u_{2\nu
}(x,z)\,dz
\end{eqnarray*}%
and this concludes the proof.

An alternative proof of the relationship (\ref{vologda}) is based on the
Fourier transforms, since for $u_{\nu }$ the following result is known:%
\[
\int_{-\infty }^{+\infty }e^{i\beta x}u_{\nu }(x,t)\,dx=E_{\nu ,1}(
-\beta ^{2}\lambda ^{2}t^{\nu }) ,
\]
where $E_{\nu ,1}( z) =\sum_{k=0}^{\infty }\frac{z^{k}}{\Gamma
(k\nu +1)}$ is the Mittag--Leffler function. Taking the Fourier transform of (%
\ref{vologda}) we get that
\begin{eqnarray*}
&&\int_{-\infty }^{+\infty }e^{i\beta x}\biggl\{ \frac{1}{\sqrt{\pi t}}%
\int_{0}^{\infty }e^{-w^{2}/(4t)}u_{2\nu }(x,w)\,dw\biggr\} \,dx \\
&&\qquad=\frac{1}{\sqrt{\pi t}}\int_{0}^{\infty }e^{-w^{2}/(4t)}E_{2\nu
,1}(-\beta ^{2}\lambda ^{2}w^{2\nu })\,dw \\
&&\qquad=\sum_{k=0}^{\infty }\frac{(-\beta ^{2}\lambda ^{2})^{k}}{\Gamma (2k\nu +1)%
}\int_{0}^{\infty }\frac{e^{-w^{2}/(4t)}}{\sqrt{\pi t}}w^{2k\nu }\,dw \\
&&\qquad=\bigl[ \mbox{for }w=2\sqrt{tz}\bigr]  \\
&&\qquad=\sum_{k=0}^{\infty }\frac{(-\beta ^{2}\lambda ^{2})^{k}}{\Gamma (2k\nu +1)%
}\frac{(2\sqrt{t})^{2k\nu +1}}{2\sqrt{\pi t}}\Gamma \biggl( \nu k+\frac{1}{2}%
\biggr)  \\
&&\qquad=\sum_{k=0}^{\infty }\frac{(-\beta ^{2}\lambda ^{2})^{k}}{\Gamma (2k\nu +1)%
}\frac{(2\sqrt{t})^{2k\nu +1}}{2\sqrt{\pi t}}\sqrt{\pi }2^{1-2\nu k}\frac{%
\Gamma ( 2\nu k) }{\Gamma ( \nu k) } \\
&&\qquad=\sum_{k=0}^{\infty }\frac{(-\beta ^{2}\lambda ^{2}t^{\nu })^{k}}{\Gamma
(k\nu +1)}=\int_{-\infty }^{+\infty }e^{i\beta x}u_{\nu }(x,t)\,dx.
\end{eqnarray*}%
\upqed
\end{pf}

\begin{remark}\label{rem2.1}
In the special case where $\nu =\frac{1}{2}$, formula (\ref{vologda}) yields
\begin{eqnarray}\label{Omsk}
u_{1/2}(x,t) &=&\frac{1}{\sqrt{\pi t}}\int_{0}^{\infty }e^{-
z^{2}/(4t)}\frac{e^{-x^{2}/(4\lambda ^{2}z)}}{\sqrt{4\pi \lambda ^{2}z}%
}\,dz \nonumber \\
&=&[2\lambda ^{2}z=y]   \\
&=&\frac{1}{\sqrt{\pi t}}\int_{0}^{\infty }\frac{e^{-x^{2}/(2y)}}{%
\sqrt{2\pi y}}\frac{e^{-y^{2}/(4t(2\lambda ^{2})^{2})}}{2\lambda ^{2}%
}\,dy.  \nonumber
\end{eqnarray}

Particularly interesting is the case where $2(2\lambda ^{2})^{2}=1,$ that is,
when $\lambda ^{2}=2^{-3/2}$, because (\ref{Omsk}) reduces to
\begin{equation}
u_{1/2}(x,t)=2\int_{0}^{\infty }\frac{e^{-x^{2}/(2y)}}{\sqrt{%
2\pi y}}\frac{e^{-y^{2}/(2t)}}{\sqrt{2\pi t}}\,dy,  \label{mezzo}
\end{equation}%
which permits us to conclude that, in this case, the solution coincides with
the probability density of the iterated Brownian motion (\ref{iterated}).
\end{remark}

\begin{remark}\label{rem2.2}
If we generalize our analysis to the $n$-dimensional case and take $\nu =%
\frac{1}{2},$ we can show that the process related to a fractional equation
of the form%
\begin{equation}\label{dim}
\frac{\partial ^{1/2}u}{\partial t^{1/2}}=\lambda ^{2}\Biggl\{ \sum_{k=1}^{n}
\frac{\partial ^{2}u}{\partial x_{k}^{2}}\Biggr\} ,\qquad x_{k}\in
\mathbb{R},t>0,
\end{equation}%
with initial condition%
\[
u_{1/2}(x_{1},x_{2},\ldots,x_{n},0)=\prod_{k=1}^{n}\delta (x_{k}),
\]
has components represented by iterated Brownian motions with a common random
time. In other words, the solution to (\ref{dim}) coincides with the
distribution of the vector process%
\[
\cases{
B_{1}(|B(t)|), \cr
\cdots \cr
B_{n}(|B(t)|),}\qquad t>0,
\]
where $B_{k},$ $k=1,\ldots,n$, are mutually independent Brownian motions and
also independent from $B$.

To check this result we evaluate the Fourier transform of the solution to (%
\ref{dim}) as follows:%
\begin{eqnarray} \label{dim.1}
&&\int_{-\infty }^{+\infty }\cdots\int_{-\infty }^{+\infty }e^{i\beta
_{1}x_{1}+\cdots+i\beta _{n}x_{n}}u_{1/2
}(x_{1},\ldots,x_{n},t)\,dx_{1}\cdots dx_{n} \nonumber \\
&&\qquad=E_{1/2,1}\Biggl( -\lambda ^{2}t^{1/2}\Biggl(
\sum_{k=1}^{n}\beta _{k}^{2}\Biggr) \Biggr)   \\
&&\qquad=\frac{2}{\sqrt{\pi }}\int_{0}^{\infty }e^{-y^{2}-2y\lambda ^{2}t^{1/2}( \sum_{k=1}^{n}\beta _{k}^{2}) }\,dy.  \nonumber
\end{eqnarray}

From (\ref{dim.1}) we get the inverse Fourier transform in the following form:%
\begin{eqnarray*}
u_{1/2}(x_{1},\ldots,x_{n},t)
&=&\frac{2}{\sqrt{\pi }}\int_{0}^{\infty }e^{-y^{2}}\prod_{k=1}^{n}\frac{e^{-
x_{k}^{2}/(2(4t^{1/2}\lambda ^{2}y))}}{\sqrt{2\pi (4t^{1/2}\lambda
^{2}y)}}\,dy \\
&=&2\int_{0}^{\infty }\frac{e^{-w^{2}/(2(2^{3}t\lambda ^{4}))}}{\sqrt{%
2\pi (2^{3}t\lambda ^{4})}}\prod_{k=1}^{n}\frac{e^{-x_{k}^{2}/(2w)}}{%
\sqrt{2\pi w}}\,dw.
\end{eqnarray*}

The main difference with respect to the case of the usual multivariate heat
equation is that the components of the iterated Brownian motions are no
longer independent because they are related to each other by the common random
time $B$ (with infinitesimal variance $2^{3}\lambda ^{4}t).$
\end{remark}

We pass now to our second theorem, which is related to the case $\nu =\frac{1%
}{2^{n}},$ $n\in \mathbb{N}$.

\begin{thm}\label{thm2.2}
For $\nu =\frac{1}{2^{n}},\lambda =2^{%
1/2^{(n+1)}-1}$ the solution to equation (\ref{equation1})
under the initial condition (\ref{condition 1}) can be written as
\begin{equation}
\label{due.6}
\qquad u_{{1/2^{n}}}(x,t)=2^{n}\int_{0}^{\infty }\cdots
\int_{0}^{\infty }\frac{e^{-x^{2}/(2z_{1})}}{\sqrt{2\pi z_{1}}}\frac{%
e^{-z_{1}^{2}/(2z_{2})}}{\sqrt{2\pi z_{2}}}\cdots \frac{e^{-%
z_{n}^{2}/(2t)}}{\sqrt{2\pi t}}\,dz_{1}\cdots dz_{n}.
\end{equation}
\end{thm}

\begin{pf}
In view of the duplication formula for the Gamma
function we can write
\begin{eqnarray}  \label{dupli}
&&\Gamma \biggl( 1-\frac{k}{2^{n+1}}-\frac{1}{2^{n+1}}\biggr) \nonumber\\[-8pt]\\[-8pt]
&&\qquad=\sqrt{\pi }2^{%
1/2^{n}+k/2^{n}}\frac{\Gamma ( 1-k/2^{n}-1/2^{n}) }{\Gamma (
1/2-k/2^{n+1}-1/2^{n+1}) }\nonumber
\end{eqnarray}%
so that the first member of (\ref{equation2}) becomes, for $\nu =\frac{1}{%
2^{n}}$ and $\lambda =2^{1/2^{(n+1)}-1}$,
\begin{eqnarray}\label{ujasnoe otnoshenie}
&&u_{1/2^{n}}(x,t)\nonumber\\
&&\qquad=\frac{1}{(2t)^{1/2^{n+1}}}\sum_{k=0}^{\infty }\biggl( -\frac{2|x|}{%
(2t)^{1/2^{n+1}}}\biggr) ^{k}\frac{1}{k!\Gamma (1-k/2^{n+1}-
1/2^{n+1})}  \nonumber \\
&&\qquad=\frac{1}{(2t)^{1/2^{n+1}}}\sum_{k=0}^{\infty }\biggl( -\frac{2|x|}{%
(2t)^{1/2^{n+1}}}\biggr) ^{k}\frac{\int_{0}^{\infty
}e^{-w_{1}}w_{1}^{-1/2^{n+1}-k/2^{n+1}-1/2}\,dw_{1}}{k!%
\sqrt{\pi }2^{(k+1)/2^{n}}\Gamma (1-k/2^{n}-1/2^{n})}
\nonumber\\
&&\qquad=\frac{1}{(2t)^{1/2^{n+1}}}\nonumber\\
&&\quad\qquad{}\times\sum_{k=0}^{\infty }\biggl( -\frac{2|x|}{%
(2t)^{1/2^{n+1}}}\biggr) ^{k}\\
&&\quad\qquad{}\times\frac{\int_{0}^{\infty
}e^{-w_{1}}w_{1}^{-1/2^{n+1}-k/2^{n+1}-1/2
}\,dw_{1}\int_{0}^{\infty }e^{-w_{2}}w_{2}^{-1/2^{n}-k/2^{n}-
1/2}\,dw_{2}}{k!(\sqrt{\pi })^{2}2^{(k+1)/2^{n}+(k+1)/
2^{n-1}}\Gamma (1-k/2^{n-1}-1/2^{n-1})}  \nonumber \\
&&\qquad=\frac{1}{(2t)^{1/2^{n+1}}}\nonumber\\
&&\quad\qquad{}\times\sum_{k=0}^{\infty }\biggl( -\frac{2|x|}{%
(2t)^{1/2^{n+1}}}\biggr) ^{k}\nonumber\\
&&\quad\qquad\phantom{{}\times\sum_{k=0}^{\infty }}
{}\times\frac{\int_{0}^{\infty
}\int_{0}^{\infty }\cdots \int_{0}^{\infty
}e^{-\sum_{j=1}^{n}w_{j}}\prod_{j=1}^{n}w_{j}^{-(k+1)/2^{n+2-j}-1/2}\,dw_{j}}{k!(\sqrt{\pi })^{n}2^{(k+1)\sum_{j=0}^{n-1}1/2^{n-j}
}\Gamma ( 1/2-k/2) }.  \nonumber
\end{eqnarray}%
At this point we can use the reflection formula for the Gamma function
\begin{eqnarray}\label{rifl}
\Gamma \biggl( \frac{1}{2}-\frac{k}{2}\biggr) &=&\frac{\pi }{\sin \{
\pi /2(1-k)\} }\frac{1}{\Gamma ( (1+k)/2) }\nonumber\\[-8pt]\\[-8pt]
&=&\frac{\pi }{\cos k\pi /2\Gamma ( (1+k)/2) }\nonumber
\end{eqnarray}%
and this shows that only even terms of (\ref{ujasnoe otnoshenie}) must\vadjust{\goodbreak} be
retained. We can therefore write that
\begin{eqnarray}\label{due.8}
&&u_{1/2^{n}}(x,t)\nonumber\\
&&\qquad=\frac{1}{(2t)^{1/2^{n+1}}}\nonumber\\
&&\qquad\quad{}\times\sum_{k=0}^{\infty }\biggl( -\frac{2|x|}{%
(2t)^{1/2^{n+1}}}\biggr) ^{k}\nonumber\\
&&\quad\qquad
{}\times\Biggl( \int_{0}^{\infty }\cdots
\int_{0}^{\infty }e^{-\sum_{j=1}^{n}w_{j}}\prod_{j=1}^{n}w_{j}^{-%
(k+1)/2^{n+2-j}-1/2}\,dw_{j}\Biggr)\nonumber \\
&&\quad\qquad{}\times \cos \frac{k\pi}{2} \Gamma
\biggl( \frac{1+k}{2}\biggr) \bigl[k!\bigl(\sqrt{\pi }\bigr)^{n}\pi
2^{(k+1)\sum_{j=0}^{n-1}1/2^{n-j}}\bigr]^{-1}  \nonumber \\
&&\qquad=\frac{2}{(2t)^{1/2^{n+1}}}\sum_{k=0}^{\infty }\biggl( -\frac{2|x|}{%
(2t)^{1/2^{n+1}}}\biggr) ^{k}\nonumber\\
&&\quad\qquad
{}\times\frac{( \int_{0}^{\infty }\cdots
\int_{0}^{\infty }e^{-\sum_{j=1}^{n}w_{j}}\prod_{j=1}^{n}w_{j}^{-%
(k+1)/2^{n+2-j}-1/2}\,dw_{j}) \cos k\pi /2\Gamma
( k) }{k!(\sqrt{\pi })^{n+1}2^{\sum_{j=0}^{n-1}1/2^{n-j}
}2^{k\sum_{j=0}^{n}1/2^{n-j}}\Gamma ( k/2) }\\
&&\qquad=\frac{2}{(2t)^{1/2^{n+1}}}\nonumber\\
&&\quad\qquad{}\times\sum_{r=0}^{\infty }\biggl( -\frac{2|x|}{%
(2t)^{1/2^{n+1}}}\biggr) ^{2r}\nonumber\\
&&\quad\qquad
{}\times\frac{( \int_{0}^{\infty }\cdots
\int_{0}^{\infty }e^{-\sum_{j=1}^{n}w_{j}}\prod_{j=1}^{n}w_{j}^{-%
(2r+1)/2^{n+2-j}-1/2}\,dw_{j}) (-1)^{r}}{(\sqrt{\pi }%
)^{n+1}2^{\sum_{j=0}^{n}1/2^{n-j}}2^{2r\sum_{j=0}^{n}1/2^{n-j}}r!}
\nonumber\\
&&\qquad=\frac{2}{(2t)^{1/2^{n+1}}2^{2( 1-1/2^{n+1}) }(%
\sqrt{\pi })^{n+1}}\nonumber\\
&&\quad\qquad{}\times\sum_{r=0}^{\infty }\frac{(-1)^{r}}{r!}\biggl( \frac{%
x^{2}}{t^{1/2^{n}}}2^{2-1/2^{n}}\biggr) ^{r}\bigl(
2^{-2( 2-1/2^{n}) }\bigr) ^{r}\nonumber\\
&&\quad\qquad
{}\times \Biggl( \int_{0}^{\infty }\cdots \int_{0}^{\infty
}e^{-\sum_{j=1}^{n}w_{j}}\prod_{j=1}^{n}w_{j}^{-(2r+1)/2^{n+2-j}-
1/2}\,dw_{j}\Biggr) .\nonumber
\end{eqnarray}

By considering that
\[
\sum_{r=0}^{\infty }\frac{(-1)^{r}}{r!}\Biggl[ \frac{x^{2}}{2^{2}}\biggl(
\frac{2}{t}\biggr) ^{1/2^{n}}\prod_{j=1}^{n}w_{j}^{-1/
2^{n+1-j}}\Biggr] ^{r} = e^{-x^{2}/2^{2}(2/t) ^{1/2^{n}
}\prod_{j=1}^{n}w_{j}^{-1/2^{n+1-j}}},
\]
we can write (\ref{due.8}) as follows:
\begin{eqnarray*}
u_{1/2^{n}}(x,t)
&=&\frac{1}{(2t)^{1/2^{n+1}}2^{1-1/2^{n}}(\sqrt{\pi })^{n+1}}\\
&&{}\times
\int_{0}^{\infty }\cdots \int_{0}^{\infty }e^{-x^{2}/2^{2}
( 2/t) ^{1/2^{n}}\prod_{j=1}^{n}w_{j}^{-1/
2^{n+1-j}}}e^{-\sum_{j=1}^{n}w_{j}}\\
&&\phantom{{}\times
\int_{0}^{\infty }\cdots \int_{0}^{\infty }}{}\times\prod_{j=1}^{n}\bigl( w_{j}^{-1/
2^{n+2-j}-1/2}\,dw_{j}\bigr).
\end{eqnarray*}

In order to calculate the integrals let us write%
\[
2( 2^{-1}t) ^{1/2^{n}}\prod_{j=1}^{n}w_{j}^{1/
2^{n+1-j}}=z_{1}
\]
so that%
\[
w_{n}=\biggl( \frac{z_{1}2^{-1}( 2^{-1}t) ^{-1/2^{n}}}{%
\prod_{j=1}^{n-1}w_{j}^{1/2^{n+1-j}}}\biggr) ^{2}
\]
and%
\[
dw_{n}=2z_{1}\,dz_{1}\biggl( \frac{2^{-1}( 2^{-1}t) ^{-1/2^{n}
}}{\prod_{j=1}^{n-1}w_{j}^{1/2^{n+1-j}}}\biggr) ^{2}.
\]

Therefore we get
\begin{eqnarray}\label{due.9}
u_{1/2^{n}}(x,t)
&=&\frac{1}{(2t)^{1/2^{n+1}}2^{1-1/2^{n}}(\sqrt{\pi
})^{n+1}}\nonumber\\
&&{}\times\int_{0}^{\infty }\cdots \int_{0}^{\infty }e^{-x^{2}/(2z_{1})
}\prod_{j=1}^{n-1}w_{j}^{-1/2^{n+2-j}-1/2} \nonumber \\
&&\phantom{{}\times\int_{0}^{\infty }\cdots \int_{0}^{\infty }}{}\times e^{-\sum_{j=1}^{n-1}w_{j}}e^{-z_{1}^{2}/(2^{2}(
2^{-1}t) ^{1/2^{n-1}}\prod_{j=1}^{n-1}w_{j}^{1/2^{n-j}}
)}\nonumber \\
&&\phantom{{}\times\int_{0}^{\infty }\cdots \int_{0}^{\infty }}{}\times 2z_{1}\biggl( \frac{2^{-1}( 2^{-1}t) ^{-1/2^{n}}}{%
\prod_{j=1}^{n-1}w_{j}^{1/2^{n+1-j}}}\biggr) ^{2}\nonumber \\
&&\phantom{{}\times\int_{0}^{\infty }\cdots \int_{0}^{\infty }}{}\times  \biggl( \frac{z_{1}2^{-1}( 2^{-1}t) ^{-1/2^{n}}}{%
\prod_{j=1}^{n-1}w_{j}^{1/2^{n+1-j}}}\biggr) ^{-1-1/2
}\,dz_{1}\,dw_{1}\cdots dw_{n-1}  \nonumber
\\
&=&\frac{\sqrt{2}( 2^{-1}t) ^{-1/2^{n+1}}}{(2t)^{1/
2^{n+1}}2^{1-1/2^{n}}(\sqrt{\pi })^{n+1}}\\
&&{}\times\int_{0}^{\infty }\frac{%
e^{-x^{2}/(2z_{1})}}{\sqrt{z_{1}}}\int_{0}^{\infty }\cdots
\int_{0}^{\infty }e^{-z_{1}^{2}/(2^{2}( 2^{-1}t) ^{1/
2^{n-1}}\prod_{j=1}^{n-1}w_{j}^{1/2^{n-j}})}  \nonumber \\
&&\phantom{{}\times\int_{0}^{\infty }}{}\times e^{-\sum_{j=1}^{n-1}w_{j}}\prod_{j=1}^{n-1}w_{j}^{-1/2^{n+1-j}-1/2}\,dz_{1}\,dw_{1}\cdots dw_{n-1}.  \nonumber
\end{eqnarray}

Now we make the similar substitution%
\[
2( 2^{-1}t) ^{1/2^{n-1}}\prod_{j=1}^{n-1}w_{j}^{1/
2^{n-j}}=z_{2}
\]
so that we get again%
\[
w_{n-1}=\biggl( \frac{z_{2}2^{-1}( 2^{-1}t) ^{-1/2^{n-1}}}{%
\prod_{j=1}^{n-2}w_{j}^{1/2^{n-j}}}\biggr) ^{2}
\]
and
\[
dw_{n-1}=2z_{2}\,dz_{2}\biggl( \frac{2^{-1}( 2^{-1}t) ^{-1/
2^{n-1}}}{\prod_{j=1}^{n-2}w_{j}^{1/2^{n-j}}}\biggr) ^{2}.
\]

In view of these substitutions, formula (\ref{due.9}) is transformed into%
\begin{eqnarray}\label{due.10}
&&u_{1/2^{n}}(x,t)\nonumber\\
&&\qquad=\frac{\sqrt{2}(2^{-1}t)^{-1/2^{n+1}}\sqrt{2}(2^{-1}t)^{-1/
2^{n}}}{(2t)^{1/2^{n+1}}2^{1-1/2^{n}}(\sqrt{\pi
})^{n+1}}\nonumber\\
&&\quad\qquad
{}\times\int_{0}^{\infty }\frac{e^{-x^{2}/(2z_{1})}}{\sqrt{z_{1}}}%
\,dz_{1}\int_{0}^{\infty }\frac{e^{-z_{1}^{2}/(2z_{2})}}{\sqrt{z_{2}}}%
\,dz_{2} \\
&&\quad\qquad{}\times\int_{0}^{\infty }\cdots \int_{0}^{\infty }e^{-z_{2}^{2}/
(2^{2}( 2^{-1}t) ^{1/2^{n-2}}\prod_{j=1}^{n-2}w_{j}^{
1/2^{n-j-1}})}\nonumber\\
&&\quad\qquad{}\times\prod_{j=1}^{n-2}w_{j}^{-1/2^{n-j}-1/2%
}e^{-\sum_{j=1}^{n-2}w_{j}}\,dw_{1}\cdots \,dw_{n-2}.  \nonumber
\end{eqnarray}

By similar transformations, after $(n-3)$ additional steps, we arrive at%
\begin{eqnarray*}
u_{1/2^{n}}(x,t)
&=&\frac{\sqrt{2^{n-1}}(2^{-1}t)^{-1/2^{n+1}-1/2^{n}-\cdots-
1/2^{3}}}{(2t)^{1/2^{n+1}}2^{1-1/2^{n}}(\sqrt{\pi }%
)^{n+1}}\\
&&{}\times\int_{0}^{\infty }\frac{e^{-x^{2}/(2z_{1})}}{\sqrt{z_{1}}}\,
dz_{1}\int_{0}^{\infty }\frac{e^{-z_{1}^{2}/(2z_{2})}}{\sqrt{z_{2}}}%
\,dz_{2}\cdots \\
&&{}\times\int_{0}^{\infty }e^{-z_{n-1}^{2}/(2^{2}[ ( 2^{-1}t)
^{1/2}w_{1}^{1/2}] )}e^{-w_{1}}w_{1}^{-1/2^{2}-
1/2}\,dw_{1}.
\end{eqnarray*}

By means of the position%
\[
2( 2^{-1}t) ^{1/2}w_{1}^{1/2}=z_{n}
\]
we get that%
\[
w_{1}=( z_{n}2^{-1}( 2^{-1}t) ^{-1/2}) ^{2}
\]
and
\[
dw_{1}=2z_{n}\,dz_{n}( 2^{-1}( 2^{-1}t) ^{-1/2})^{2}.
\]

We arrive at the final expression
\begin{eqnarray*}
u_{1/2^{n}}(x,t)
&=&\frac{\sqrt{2^{n}}(2^{-1}t)^{-1/2^{n+1}-1/2^{n}-\cdots-
1/2^{3}-1/2^{2}}}{(2t)^{1/2^{n+1}}2^{1-1/2^{n}}(%
\sqrt{\pi })^{n+1}} \\
&&{}\times\int_{0}^{\infty }\frac{e^{-x^{2}/(2z_{1})}}{\sqrt{z_{1}}}%
\,dz_{1}\int_{0}^{\infty }\frac{e^{-z_{1}^{2}/(2z_{2})}}{\sqrt{z_{2}}}%
\,dz_{2}\cdots\\
&&{}\times \int_{0}^{\infty }\frac{e^{-z_{n-1}^{2}/(2z_{n})%
}}{\sqrt{z_{n}}}e^{-z_{n}^{2}/(2t)}\,dz_{n} \\
&=&\frac{2^{n}}{2^{n/2+1/2}(\sqrt{\pi })^{n+1}\sqrt{t}}\\
&&{}\times
\int_{0}^{\infty
}\frac{e^{-x^{2}/(2z_{1})}}{\sqrt{z_{1}}}\,dz_{1}\cdots
\int_{0}^{\infty }\frac{e^{-z_{n-1}^{2}/(2z_{n})}}{\sqrt{%
z_{n}}}e^{-z_{n}^{2}/(2t)}\,dz_{n},
\end{eqnarray*}%
which coincides with (\ref{due.6}).
\end{pf}

\begin{remark}\label{rem2.3}
It is well known that the Laplace--Fourier transform of the solution to (\ref%
{equation1}) with initial conditions (\ref{condition 1}) or (\ref{cond1}) is
equal, for $0<\nu \leq 2,$ to%
\begin{equation}\label{lf}
\qquad\int_{0}^{+\infty }e^{-st}\,ds\int_{-\infty }^{+\infty }e^{i\beta x}u_{\nu
}(x,t)\,dx=\frac{s^{\nu -1}}{s^{\nu }+\lambda ^{2}\beta ^{2}},\qquad s>0,\beta
\in \mathbb{R}.
\end{equation}

We check that the Laplace--Fourier transform of (\ref{due.6}) reduces to (\ref%
{lf}) for $\nu =\frac{1}{2^{n}}$ and $\lambda ^{2}=2^{1/2^{n}-2}$:
\begin{eqnarray*}
&&\int_{-\infty }^{+\infty }e^{i\beta x}u_{1/2^{n}
}(x,t)\,dx\\
&&\qquad=2^{n}\int_{-\infty }^{+\infty }e^{i\beta x}\,dx\int_{0}^{\infty }%
\frac{e^{-x^{2}/(2z_{1})}}{\sqrt{2\pi z_{1}}}\,dz_{1}\cdots
\int_{0}^{\infty }\frac{e^{-z_{n}^{2}/(2t)}}{\sqrt{2\pi t}}\,dz_{n} \\
&&\qquad=2^{n}\int_{0}^{\infty }e^{-\beta ^{2}/2z_{1}}\frac{e^{-
z_{1}^{2}/(2z_{2})}}{\sqrt{2\pi z_{2}}}\,dz_{1}\int_{0}^{\infty }\frac{e^{-%
z_{2}^{2}/(2z_{3})}}{\sqrt{2\pi z_{3}}}\,dz_{2}\cdots
\int_{0}^{\infty }\frac{e^{-z_{n}^{2}/(2t)}}{\sqrt{2\pi t}}\,dz_{n}
\\
&&\qquad=2^{n}\sum_{r=0}^{\infty }\biggl( -\frac{\beta ^{2}}{2}\biggr) ^{r}\frac{1}{%
r!}\int_{0}^{\infty }z_{1}^{r}\frac{e^{-z_{1}^{2}/(2z_{2})}}{\sqrt{%
2\pi z_{2}}}\,dz_{1}\cdots \int_{0}^{\infty }\frac{e^{-
z_{n}^{2}/(2t)}}{\sqrt{2\pi t}}\,dz_{n} \\
&&\qquad=2^{n}\sum_{r=0}^{\infty }\biggl( -\frac{\beta ^{2}}{2}\biggr) ^{r}\frac{1}{%
r!}\frac{2^{r/2-1}}{\sqrt{\pi }}\Gamma \biggl( \frac{r+1}{2}\biggr)\\
&&\quad\qquad\phantom{2^{n}\sum_{r=0}^{\infty }}
{}\times
\int_{0}^{\infty }z_{2}^{r/2}\frac{e^{-z_{2}^{2}/(2z_{3})}}{%
\sqrt{2\pi z_{3}}}\,dz_{2}\cdots \int_{0}^{\infty }\frac{e^{-
z_{n}^{2}/(2t)}}{\sqrt{2\pi t}}\,dz_{n} \\
&&\qquad=2^{n}\sum_{r=0}^{\infty }\biggl( -\frac{\beta ^{2}}{2}\biggr) ^{r}\frac{1}{%
r!}\frac{2^{r/2-1}2^{r/4-1}}{( \sqrt{\pi }) ^{2}}%
\Gamma \biggl( \frac{r}{2}+\frac{1}{2}\biggr) \Gamma \biggl( \frac{r}{4}+\frac{%
1}{2}\biggr) \\
&&\qquad\quad\phantom{2^{n}\sum_{r=0}^{\infty }}
{}\times\int_{0}^{\infty }z_{3}^{r/4}\frac{e^{-z_{3}^{2}/
(2z_{4})}}{\sqrt{2\pi z_{4}}}\,dz_{3}\cdots \int_{0}^{\infty }\frac{%
e^{-z_{n}^{2}/(2t)}}{\sqrt{2\pi t}}\,dz_{n}
\\
&&\qquad=2^{n}\sum_{r=0}^{\infty }\biggl( -\frac{\beta ^{2}}{2}\biggr) ^{r}\frac{1}{%
r!}\frac{2^{r/2-1}2^{r/4-1}\cdots 2^{r/2^{n-1}-1}}{%
( \sqrt{\pi }) ^{n-1}}\\
&&\qquad\quad\phantom{2^{n}\sum_{r=0}^{\infty }}{}\times\Gamma \biggl( \frac{r}{2}+\frac{1}{2}\biggr) \Gamma \biggl( \frac{r}{%
2^{2}}+\frac{1}{2}\biggr) \cdots \Gamma \biggl( \frac{r}{2^{n-1}}+\frac{1%
}{2}\biggr)\\
&&\qquad\quad\phantom{2^{n}\sum_{r=0}^{\infty }}{}\times \int_{0}^{\infty }z_{n}^{r/2^{n-1}}\frac{e^{-
z_{n}^{2}/(2t)}}{\sqrt{2\pi t}}\,dz_{n} \\
&&\qquad=2^{n}\sum_{r=0}^{\infty }\biggl( -\frac{\beta ^{2}}{2}\biggr) ^{r}\frac{1}{%
r!}\frac{2^{r/2+r/4+\cdots+r/2^{n}-n}}{( \sqrt{\pi
}) ^{n}}t^{r/2^{n}}\\
&&\quad\qquad\\
&&\qquad\quad\phantom{2^{n}\sum_{r=0}^{\infty }}{}\times \Gamma \biggl( \frac{r}{2}+\frac{1}{2}%
\biggr) \Gamma \biggl( \frac{r}{2^{2}}+\frac{1}{2}\biggr) \cdots \Gamma
\biggl( \frac{r}{2^{n}}+\frac{1}{2}\biggr) .
\end{eqnarray*}

By applying the duplication formula we get that%
\begin{eqnarray} \label{ast}
&&\Gamma \biggl( \frac{r}{2}+\frac{1}{2}\biggr) \Gamma \biggl( \frac{r}{2^{2}}+%
\frac{1}{2}\biggr) \cdots \Gamma \biggl( \frac{r}{2^{n}}+\frac{1}{2}%
\biggr)\nonumber \\
&&\qquad=\sqrt{\pi }2^{1-r}\frac{\Gamma (r)}{\Gamma ( r/2) }%
\sqrt{\pi }2^{1-r/2}\frac{\Gamma ( r/2) }{\Gamma
( r/2^{2}) }\cdots \sqrt{\pi }2^{1-r/
2^{n-1}}\frac{\Gamma ( r/2^{n-1}) }{\Gamma ( r/2^{n}) }   \\
&&\qquad=\sqrt{\pi ^{n}}2^{n-r-r/2-\cdots-r/2^{n-1}}\frac{\Gamma (r)}{%
\Gamma ( r/2^{n}) }  \nonumber
\end{eqnarray}
\noindent and thus%
\begin{eqnarray}\label{due.11}
&&\int_{-\infty }^{+\infty }e^{i\beta x}u_{1/2^{n}}(x,t)\,dx \nonumber\\
&&\qquad=2^{n}\sum_{r=0}^{\infty }\biggl( -\frac{\beta ^{2}}{2}\biggr) ^{r}\frac{1}{%
r!}2^{r/2+r/4+\cdots+r/2^{n}-n}2^{n-r-r/2-\cdots-
r/2^{n-1}}t^{r/2^{n}}\frac{\Gamma (r)}{\Gamma ( r/2^{n}) } \!\!\!\!\!\!\!\! \nonumber \\
&&\qquad=\sum_{r=0}^{\infty }\biggl( -\frac{\beta ^{2}}{2}\biggr)^{r}\frac{2^{r/2^{n}-r}t^{r/2^{n}}}{r/2^{n}\Gamma (
r/2^{n}) }   \\
&&\qquad=\sum_{r=0}^{\infty }\biggl( -\frac{\beta ^{2}t^{1/2^{n}}}{2^{2-%
1/2^{n}}}\biggr) ^{r}\frac{1}{\Gamma ( r/2^{n}+1)
}  \nonumber \\
&&\qquad=E_{1/2^{n},1}\biggl( -\frac{\beta ^{2}t^{1/2^{n}}}{2^{2-%
1/2^{n}}}\biggr) .  \nonumber
\end{eqnarray}

By taking the Laplace transform of (\ref{due.11}) we get%
\[
\int_{0}^{+\infty }e^{-st}E_{1/2^{n},1}\biggl( -\frac{\beta ^{2}t^{%
1/2^{n}}}{2^{2-1/2^{n}}}\biggr) \,dt=\frac{s^{1/2^{n}%
-1}2^{2-1/2^{n}}}{\beta ^{2}+2^{2-1/2^{n}}s^{1/2^{n}}},
\]
which coincides with (\ref{lf}), for $\nu =\frac{1}{2^{n}}$ and $\lambda
^{2}=2^{1/2^{n}-2}$.

The form (\ref{due.6}) of the solution $u_{1/2^{n}}$ shows that it
coincides with the distribution of the $n$-times iterated Brownian motion
defined in (\ref{iperiterated}).
\end{remark}

Another representation of the solution to the fractional equation (\ref%
{equation1}) can be inferred from the following result:

\begin{thm}\label{thm2.3}
The solution $u_{\nu }(x,t)=u_{\nu }$
to the initial value problem (\ref{eq4}), for $0<\nu \leq 1$,
can be written as
\begin{equation} \label{new.1}
u_{\nu }(x,t)=\int_{0}^{\infty }\frac{1}{\sqrt{4\pi w\lambda }}e^{-
x^{2}/(4w\lambda) }\overline{u}_{2\nu }(w,t)\,dw,
\end{equation}%
where
\begin{equation}
\label{new.2}
\overline{u}_{2\nu }(w,t)=%
\cases
{
2u_{2\nu }(w,t),&\quad for $w\geq 0$, \cr
0,&\quad for $w<0$
}
\end{equation}
and $u_{2\nu }$ is the solution of (\ref{init1}) or (\ref%
{init2}).
\end{thm}

\begin{pf}
We first note that for the solutions to (\ref{init1}%
) or (\ref{init2}) the following result holds:%
\begin{equation}\label{lt}
L(x,s) =\int_{0}^{\infty }e^{-st}u_{2\nu }(x,t)\,dt
=\frac{s^{\nu -1}}{2\lambda }e^{-|x|s^{\nu }/\lambda },
\end{equation}%
as can be obtained by taking the Laplace transform of $\frac{\partial ^{2\nu
}u}{\partial t^{2\nu }}=\lambda ^{2}\frac{\partial ^{2}u}{\partial x^{2}}.$
The solution to the corresponding equation
\[
s^{2\nu }L-s^{2\nu -1}\delta (x)=\lambda ^{2}\frac{d^{2}L}{dx^{2}}
\]
coincides with the solution to%
\[
\cases{
\lambda ^{2}\dfrac{d^{2}L}{dx^{2}}=s^{2\nu }L,&\quad $x\neq 0$, \cr
\dfrac{dL}{dx}\Bigm\vert ^{+}- \dfrac{dL}{dx}\Bigm\vert ^{-}=-%
\dfrac{s^{2\nu -1}}{\lambda ^{2}},\cr
L(s,0^{+})=L(s,0^{-}),
}
\]
and easily yields (\ref{lt}); see also (3.3) of Orsingher and Beghin
(\citeyear{OrsBeg2004}). Therefore, by taking the Laplace transform of (\ref{new.1}), we get%
\begin{eqnarray*}
&&\int_{0}^{\infty }\frac{1}{\sqrt{4\pi w\lambda }}e^{-x^{2}/
(4w\lambda)}\biggl\{ 2\int_{0}^{\infty }e^{-st}u_{2\nu }(w,t)\,dt\biggr\} \,dw \\
&&\qquad=2\int_{0}^{\infty }\frac{1}{\sqrt{4\pi w\lambda }}e^{-x^{2}/
(4w\lambda)}\frac{s^{\nu -1}}{2\lambda }e^{-s^{\nu }/\lambda w} \,dw \\
&&\qquad=[ 2w=z] \\
&&\qquad=\frac{s^{\nu -1}}{2\lambda }\int_{0}^{\infty }\frac{1}{\sqrt{2\pi
z\lambda }}e^{-x^{2}/(2z\lambda)}e^{-s^{\nu }/\lambda z/2}\,dz \\
&&\qquad=\frac{s^{\nu /2-1}}{2\lambda }e^{-|x|s^{\nu /2}
/\lambda }
\end{eqnarray*}%
and this coincides with the Laplace transform of $u_{\nu }(x,t)$.
\end{pf}

\begin{remark}\label{rem2.4}
Formula (\ref{new.1}) suggests that we should represent the solution of (\ref{eq4})
as the distribution of the process%
\begin{equation} \label{new.3}
B(\mathcal{T}_{2\nu }(t)),\qquad t>0,
\end{equation}%
where $B$ is a Brownian motion with infinitesimal variance $2\lambda $ and $%
\mathcal{T}_{2\nu }(t)$, $t>0$, is a process, independent from $B,$ with law
equal to (\ref{new.2}).

It is straightforward that, for $\nu =1/2$, the process (\ref{new.3})
coincides with the iterated Brownian motion \emph{$\mathcal{I}$}$_{1}$; see (%
\ref{iterated}).

By comparing the relationship (\ref{new.1}) with (\ref{vologda}) we note
also that, in the composition of processes, Brownian motion plays in the
second case the role of ``time,'' while
in the first one it represents ``space.''
\end{remark}

\section{On moments and functionals of the iterated Brownian
motion}\label{s3}

Some properties of the classical iterated Brownian motion have been obtained
by several authors and include the law of iterated logarithm [Burdzy and San
Mart\`{\i}n (\citeyear{BurMar1995})] and the modulus of continuity [Khoshnevisan and Lewis
(\citeyear{KhoLew1996})]. Applications of the iterated Brownian motion to diffusion in cracks
are dealt with in De Blassie (\citeyear{DeB2004}).

We start by presenting the distribution of the maximum of the $n$-times
iterated Brownian motion and, in an explicit form, for the usual iterated
Brownian motion.

\begin{thm}\label{thm3.1}
For the $n$-times iterated Brownian motion
\[
\mathcal{I}_{n}(t)=B_{1}(|B_{2}(|B_{3}(\cdots(|B_{n+1}(t)|)\cdots)|)|),\qquad t>0,
\]
where $B_{j}$, $j=1,\ldots,n+1$, are independent Brownian
motions, we have for $\beta >0$ that
\begin{eqnarray}\label{max.1}
&&\Pr \biggl\{ \max_{0\leq s\leq t}\mathcal{I}_{n}(s)\in d\beta \biggr\} \nonumber\\
&&\qquad=2\int_{0}^{+\infty }\cdots\int_{0}^{+\infty }\Pr \{ B
_{1}(y_{1})\in d\beta \} \Pr \biggl\{ \max_{0\leq z_{1}\leq y_{2}}|%
B_{2}(z_{1})|\in dy_{1}\biggr\}  \\
&&\quad\qquad{}\times\Pr \biggl\{ \max_{0\leq z_{2}\leq y_{3}}|B_{3}(z_{2})|\in
dy_{2}\biggr\} \cdots \Pr \biggl\{ \max_{0\leq z_{n}\leq t}|\emph{B%
}_{n+1}(z_{n})|\in dy_{n}\biggr\} .  \nonumber
\end{eqnarray}
\end{thm}

\begin{pf}
For  $\mathcal{I}_{1}(t)=B_{1}(|B_{2}(t)|) $ we can write that%
\begin{eqnarray}\label{max.2}
&&\Pr \biggl\{ \max_{0\leq s\leq t}\mathcal{I}_{1}(s)\in d\beta \biggr\}\nonumber
\\
&&\qquad=\Pr \biggl\{ \max_{0\leq z\leq \max_{0\leq w\leq t}|B_{2}(w)|}
B_{1}(z)\in d\beta \biggr\}  \nonumber
\\
&&\qquad=E\biggl\{ \Pr \biggl\{  \max_{0\leq z\leq \max_{0\leq w\leq
t}|B_{2}(w)|}B_{1}(z)\in d\beta \Bigm\vert \max_{0\leq w\leq
t}|B_{2}(w)|\biggr\} \biggr\}   \\
&&\qquad=\int_{0}^{+\infty }\Pr \biggl\{ \max_{0\leq z\leq y}B_{1}(z)\in
d\beta \biggr\} \Pr \biggl\{ \max_{0\leq w\leq t}|B_{2}(w)|\in
dy\biggr\}  \nonumber \\
&&\qquad=2\int_{0}^{+\infty }\Pr \{ B_{1}(y)\in d\beta \} \Pr
\biggl\{ \max_{0\leq w\leq t}|B_{2}(w)|\in dy\biggr\} .  \nonumber
\end{eqnarray}

For  $\mathcal{I}_{n}(t)=B_{1}(|\mathcal{I}_{n-1}(t)|),$ $n\geq 1$, we have analogously that
\begin{eqnarray}\label{max.3}
&&\Pr \biggl\{ \max_{0\leq s\leq t}\mathcal{I}_{n}(s)\in d\beta \biggr\}
\nonumber\\[-8pt]\\[-8pt]
&&\qquad=2\int_{0}^{+\infty }\Pr \{ B_{1}(y)\in d\beta \} \Pr
\biggl\{ \max_{0\leq w\leq t}|\mathcal{I}_{n-1}(w)|\in dy\biggr\}  \nonumber
\end{eqnarray}
and, by induction, we obtain (\ref{max.1}).
\end{pf}

\begin{remark}\label{rem3.1}
In the case $n=1$ we can give an explicit expression for (\ref{max.1}) as
follows:%
\begin{eqnarray}\label{max.4}
&&\Pr \biggl\{ \max_{0\leq s\leq t}\mathcal{I}_{1}(s)\in d\beta \biggr\}\nonumber \\
&&\qquad=2d\beta \int_{0}^{+\infty }\frac{e^{-\beta ^{2}/(2w)}}{\sqrt{2\pi w}%
}\nonumber\\
&&\qquad\phantom{=2d\beta \int_{0}^{+\infty }}
{}\times\Biggl\{ \sum_{k=-\infty }^{+\infty }(-1)^{k}\biggl[ (1+2k)\frac{e^{-
w^{2}/(2t)(1+2k)^{2}}}{\sqrt{2\pi t}}\nonumber\\[-8pt]\\[-8pt]
&&\qquad\quad\hspace*{43pt}\phantom{{}\times\Biggl\{ \sum_{k=-\infty }^{+\infty }(-1)^{k}\biggl[}
{}+(1-2k)\frac{e^{-w^{2}/(2t)
(1-2k)^{2}}}{\sqrt{2\pi t}}\biggr] \Biggr\} \,dw  \nonumber \\
&&\qquad=2\sum_{k=-\infty }^{+\infty }(-1)^{k}\biggl[ \Pr \biggl\{ \mathcal{I}%
_{1}\biggl( \frac{t}{(1+2k)^{2}}\biggr) \in d\beta \biggr\} +\Pr \biggl\{
\mathcal{I}_{1}\biggl( \frac{t}{(1-2k)^{2}}\biggr) \in d\beta \biggr\} \biggr]
\nonumber \\
&&\qquad=2d\beta \sum_{k=-\infty }^{+\infty }(-1)^{k}\biggl[ u_{1/2}\biggl(
\beta ,\frac{t}{(1+2k)^{2}}\biggr) +u_{1/2}\biggl( \beta ,\frac{t}{%
(1-2k)^{2}}\biggr) \biggr] ,  \nonumber
\end{eqnarray}%
where $u_{1/2}(x,t)$ is given in (\ref{mezzo}) and in the first step
we applied the well-known result for the maximal distribution of the
absolute value of Brownian motion [see Shorack and Wellner (\citeyear{ShoWel1986}), page 34].
The last term of (\ref{max.4}) shows that the distribution of the maximum of
the iterated Brownian motion can be expressed in terms of its probability
law $u_{1/2}=u_{1/2}(x,t)$, as in the case of the classical
Brownian motion.

In principle we could write explicitly the distribution of the maximum of
$\mathcal{I}_{n}(t)$ in terms of $u_{1/2^{n}}$, but this
produces a sum of $2^{n}$ terms, each of which has a very entangled
structure.
\end{remark}

On the basis of the same principles it is possible to write down the
distribution of the sojourn time on the positive half-line of the process $%
\mathcal{I}_{n}(t)=B_{1}(|\mathcal{I}_{n-1}(t)|),$ $t>0,$ $n\geq 1,$ defined
as
\begin{equation}\label{soj.1}
\Gamma _{t}=\int_{0}^{\max_{0\leq w\leq t}|\mathcal{I}%
_{n-1}(w)|}1_{\{ z:B_{1}(z)>0\} }\,dz.
\end{equation}

This random variable takes values in $[ 0,+\infty ) $, because
during the interval $[ 0,t) $ the process $|\mathcal{I}_{n-1}|$
(which plays the role of time for $B_{1}$) can span the whole positive real
axes.

\begin{thm}
For the process $\mathcal{I}_{n}(t)$, $t>0$, the distribution of
$\Gamma _{t}$ reads
\begin{eqnarray}\label{sogg}
\Pr \{ \Gamma _{t}\in ds\} =ds\int_{s}^{+\infty }\frac{1}{\pi
\sqrt{s(z-s)}}\Pr \biggl\{ \max_{0\leq w\leq t}|\mathcal{I}%
_{n-1}(w)|\in dz\biggr\} ,\nonumber\\[-10pt]\\[-10pt]
\eqntext{0\leq s<\infty .}
\end{eqnarray}
\end{thm}

\begin{pf}
The definition of $\Gamma _{t}$ given in (\ref{soj.1}) implies that
\begin{eqnarray} \label{soj.2}
&&\Pr \{ \Gamma _{t}\in ds\} \nonumber\\
&&\qquad=E\Biggl\{ \Pr \Biggl\{  \Biggl[ \int_{0}^{\max_{0\leq w\leq t}|%
\mathcal{I}_{n-1}(w)|}1_{\{ z:B_{1}(z)>0\} }\,dz\Biggr] \in
ds\biggm\vert \max_{0\leq w\leq t}|\mathcal{I}_{n-1}(w)|\Biggr\}
\Biggr\}  \!\!\!\!\!\!\!\! \\
&&\qquad=\int_{s}^{+\infty }\Pr \{ \Gamma _{z}\in ds\} \Pr \biggl\{
\max_{0\leq w\leq t}|\mathcal{I}_{n-1}(w)|\in dz\biggr\} .  \nonumber
\end{eqnarray}

By inserting the arc-sine law in (\ref{soj.2}) we get (\ref{sogg}).

We can check that (\ref{sogg}) integrates to one%
\begin{eqnarray*}
&&\int_{0}^{+\infty }\Pr \{ \Gamma _{t}\in ds\} \\
&&\qquad=\int_{0}^{+\infty }ds\int_{s}^{+\infty }\frac{1}{\pi \sqrt{s(z-s)}}\Pr
\biggl\{ \max_{0\leq w\leq t}|\mathcal{I}_{n-1}(w)|\in dz\biggr\} \\
&&\qquad=\int_{0}^{+\infty }\Pr \biggl\{ \max_{0\leq w\leq t}|\mathcal{I}%
_{n-1}(w)|\in dz\biggr\} \int_{0}^{z}\frac{ds}{\pi \sqrt{s(z-s)}}=1.
\end{eqnarray*}
\upqed
\end{pf}

\begin{remark}\label{rem3.2}
For the iterated Brownian motion $\mathcal{I}_{1}(t)=B_{1}(|B_{2}(t)|)$ the
distribution of $\Gamma _{t}$ can be written explicitly as follows:%
\begin{eqnarray} \label{soj.3}
&&\Pr \{ \Gamma _{t}\in ds\} \nonumber\\
&&\qquad=ds\int_{s}^{+\infty }\frac{dz}{\pi \sqrt{s(z-s)}}\nonumber\\
&&\quad\qquad
{}\times\sum_{k=-\infty
}^{+\infty }(-1)^{k}\biggl\{ \frac{e^{-z^{2}/(2t)(1+2k)^{2}}}{\sqrt{%
2\pi t}}(1+2k)+\frac{e^{-z^{2}/(2t)(1-2k)^{2}}}{\sqrt{2\pi t}}%
(1-2k)\biggr\}   \\
&&\qquad=\frac{ds}{\pi \sqrt{2\pi ts}}\sum_{k=-\infty }^{+\infty }(-1)^{k}\biggl\{
(1+2k)\int_{s}^{+\infty }\frac{e^{-z^{2}/(2t)(1+2k)^{2}}}{\sqrt{z-s}}%
\,dz\nonumber\\
&&\hspace*{135pt}
{}+(1-2k)\int_{s}^{+\infty }\frac{e^{-z^{2}/(2t)(1-2k)^{2}}}{\sqrt{z-s}%
}\,dz\biggr\} .  \nonumber
\end{eqnarray}

By the transformation $z=s(1+x^{2})$ the integrals in (\ref{soj.3}) are
converted [for $A=\frac{(1\pm 2k)^{2}}{2t}$] into%
\begin{eqnarray*}
&&2s\int_{0}^{+\infty }\frac{e^{-s^{2}A(1+x^{2})^{2}}}{\sqrt{s}}\,dx \\
&&\qquad=2\sqrt{s}e^{-s^{2}A}\int_{0}^{+\infty }e^{-s^{2}A(x^{4}+2x^{2})}\,dx \\
&&\qquad=\sqrt{\frac{s}{2}}e^{-s^{2}A/2}K_{1/4}\biggl( \frac{As^{2}%
}{2}\biggr) ,
\end{eqnarray*}%
where, in the last step, we have applied formula 3.469.1 of Gradshteyn and
Rhyzik (\citeyear{GraRyz1994}) and $K_{1/4}(x)=\frac{\pi }{\sqrt{2}}[ I_{-1/4}(x)-I_{1/4}(x)] $ [by formula 8.485 of Gradshteyn and
Rhyzik (\citeyear{GraRyz1994})]. By $I_{\nu }$ we denote the Bessel function of imaginary
argument of order $\nu $, that is, $I_{\nu }(x)=\sum_{k=0}^{+\infty }\frac{%
(x/2)^{2k+\nu }}{k!\Gamma (k+\nu +1)}$.\vadjust{\goodbreak} Therefore we get
\begin{eqnarray*}
&&\Pr \{ \Gamma _{t}\in ds\} \\
&&\qquad=\frac{ds}{2\pi \sqrt{\pi t}}\sum_{k=-\infty }^{+\infty }(-1)^{k}\biggl\{%
(1+2k)e^{-s^{2}/(4t)(1+2k)^{2}}K_{1/4}\biggl( \frac{%
s^{2}(1+2k)^{2}}{4t}\biggr) \\
&&\hspace*{95pt}\quad\qquad{}+(1-2k)e^{-s^{2}/(4t)(1-2k)^{2}}K_{1/4}\biggl( \frac{%
s^{2}(1-2k)^{2}}{4t}\biggr) \biggr\}.
\end{eqnarray*}

We now derive the explicit form of the moments of even order of \emph{$%
\mathcal{I}$}$_{n}(t)$.
\end{remark}

\begin{thm}\label{thm3.3}
For the process $\mathcal{I}_{n}(t)$, $t>0$, the moments of order
$2k$ are given by
\begin{eqnarray}\label{due.14}
E\mathcal{I}_{n}^{2k}(t)
&=&\frac{(2k)!}{k!}\frac{2^{n}}{2^{k}}\int_{0}^{\infty
}x^{k}\,dx\int_{0}^{\infty }\frac{e^{-x^{2}/(2z_{1})}}{\sqrt{2\pi z_{1}}}%
\,dz_{1}\cdots \int_{0}^{\infty }\frac{e^{-z_{n-1}^{2}/(2t)}}{%
\sqrt{2\pi t}}\,dz_{n-1}  \nonumber\\[-8pt]\\[-8pt]
&=&\frac{2^{k/2^{n}}}{2^{2k}}\frac{(2k)!}{\Gamma ( k/2^{n}+1) }t^{k/2^{n}}.  \nonumber
\end{eqnarray}
\end{thm}

\begin{pf}
The first expression in (\ref{due.14}) can be
proved by observing that, for \mbox{$n\geq 1$},
\begin{eqnarray}\label{due.15}
E\mathcal{I}_{n}^{2k}(t)
&=&E[ B_{1}^{2k}(|B_{2}(|B_{3}(\cdots|B_{n+1}(t)|\cdots)|)|)]  \nonumber \\
&=&\frac{(2k)!}{k!}\frac{1}{2^{k}}E\vert
B_{2}(|B_{3}(\cdots|B_{n+1}(t)|\cdots)|)\vert ^{k}
\nonumber\\[-8pt]\\[-8pt]
&=&\frac{(2k)!}{k!}\frac{1}{2^{k}}2\int_{0}^{+\infty }x^{k}\Pr \{
B_{2}(|B_{3}(\cdots|B_{n+1}(t)|\cdots)|)\in dx\}  \nonumber \\
&=&\frac{(2k)!}{k!}\frac{1}{2^{k}}2\int_{0}^{+\infty }x^{k}\Pr \{
\mathcal{I}_{n-1}(t)\in dx\} ,  \nonumber
\end{eqnarray}%
which coincides with the second line of (\ref{due.14}). By performing the
integrations in (\ref{due.15}) we get the explicit expression of the moments
of order $2k$:%
\begin{eqnarray}\label{due.16}
E\mathcal{I}_{n}^{2k}(t)
&=&\frac{\Gamma ( k/2+1/2) \Gamma ( k/2^{2}+1/2) \cdots \Gamma (
k/2^{n}+1/2) 2^{k/2+\cdots+k/2^{n}}t^{k/2^{n}}}{2^{n}\sqrt{%
\pi ^{n}}}\nonumber\\
&&{}\times\frac{(2k)!2^{n}}{k!2^{k}}  \nonumber \\
&=&[ \mbox{by (\ref{ast})}]   \\
&=&\sqrt{\pi ^{n}}2^{n-k-k/2-\cdots-k/2^{n-1}}\frac{\Gamma (k)}{%
\Gamma ( k/2^{n}) }\frac{2^{k/2+\cdots+k/2^{n}%
}t^{k/2^{n}}}{2^{k}\sqrt{\pi ^{n}}}\frac{(2k)!}{k!}  \nonumber \\
&=&t^{k/2^{n}}2^{n-2k+k/2^{n}}\frac{(2k)!}{k\Gamma (
k/2^{n}) }.  \nonumber
\end{eqnarray}%
\upqed
\end{pf}

\begin{remark}\label{rem3.3}
For $n=0$ formula (\ref{due.14}) coincides with the moments $EB^{2k}(t)$, which is as
it should be, since \emph{$\mathcal{I}$}$_{0}(t)=B(t).$

For $n=1$, the moments of the iterated Brownian motion \emph{$\mathcal{I}$}$%
_{1}(t)=B_{1}(|B_{2}(t)|)$ can be evaluated directly as follows:%
\begin{eqnarray*}
E\mathcal{I}_{1}^{2k}(t)
&=&EB_{1}^{2k}(|B_{2}(t)|) \\
&=&\frac{(2k)!}{k!}\frac{1}{2^{k}}E\vert B_{2}(t)\vert ^{k} \\
&=&\frac{(2k)!}{k!}\frac{2}{2^{k}}\int_{0}^{+\infty }x^{k}\frac{e^{-
x^{2}/(2t)}}{\sqrt{2\pi t}}\,dx \\
&=&\frac{2^{k/2}}{2^{2k}}\frac{(2k)!}{\Gamma ( k/2
+1) }t^{k/2},
\end{eqnarray*}%
which coincides with (\ref{due.14}) for $n=1.$

For any $n\geq 1$ and $k=1$, we obtain the explicit form of the variance%
\[
\operatorname{var}\mathcal{I}_{n}(t)=\frac{2^{1/2^{n}}t^{1/2^{n}}}{2\Gamma
( 1/2^{n}+1) },
\]
while, for $n=0,$ it is $\operatorname{var}\mathcal{I}_{0}(t)=t$, as expected.
\end{remark}

\begin{remark}\label{rem3.4}
For all $t>0$, the sequence \emph{$\mathcal{I}$}$_{n}(t)$ converges in
distribution$,$ for $n\rightarrow \infty $, to the Gauss--Laplace exponential
random variable and its density is independent from $t.$ From (\ref{multiple
integral bis}) we get that%
\begin{equation}
\lim_{n\rightarrow \infty }u_{1/2^{n}}(x,t)=e^{-2|x|},\qquad t>0,x\in
\mathbb{R}.  \label{lim.1}
\end{equation}

By working on the Fourier transform (\ref{due.11}) of $u_{1/2^{n}}$
we have the following alternative proof:%
\begin{eqnarray}\label{lim.4}
\qquad\quad\lim_{n\rightarrow \infty }\int_{-\infty }^{+\infty }e^{i\beta x}u_{1/2^{n}}(x,t)\,dx
&=&E_{0,1}\biggl( -\frac{\beta ^{2}}{2^{2}}\biggr)
=\sum_{k=0}^{\infty }\biggl( -\frac{\beta ^{2}}{2^{2}}\biggr) ^{k}=\frac{%
2^{2}}{2^{2}+\beta ^{2}}.
\end{eqnarray}

Formula (\ref{lim.4}) coincides with the characteristic function of (\ref%
{lim.1}). Loosely speaking, this shows that the composition of infinite
Brownian motions produces the bilateral exponential distribution.

In view of (\ref{multiple integral}) we have also the identity%
\begin{eqnarray}\label{lim.2}
&&\lim_{n\rightarrow \infty }2^{n}\int_{0}^{\infty }\cdots
\int_{0}^{\infty }\frac{e^{-x^{2}/(2z_{1})}}{\sqrt{2\pi z_{1}}}\frac{%
e^{-z_{1}^{2}/(2z_{2})}}{\sqrt{2\pi z_{2}}}\cdots \frac{e^{-%
z_{n}^{2}/(2t)}}{\sqrt{2\pi t}}\,dz_{1}\cdots
dz_{n}\nonumber\\[-8pt]\\[-8pt]
&&\qquad=e^{-2|x|},\nonumber
\end{eqnarray}
which is a rather striking result. Furthermore, if we assume that%
\[
\lim_{n\rightarrow \infty }\frac{\partial ^{1/2^{n}}u}{\partial t^{%
1/2^{n}}}=u,
\]
the fractional equation (\ref{equation1}) is converted into%
\[
u=\frac{1}{2^{2}}\frac{\partial ^{2}u}{\partial x^{2}},
\]
subject to%
\[
u(x,0)=\delta (x),
\]
which is satisfied by (\ref{lim.1}) for all $x\neq 0$.
\end{remark}

\begin{remark}\label{rem3.5}
For the random process
\begin{equation}
T( |B_{2}(|B_{3}(\cdots|B_{n+1}(t)|\cdots)|)|) ,\qquad t>0,
\label{tel.1}
\end{equation}
where $T$ is a telegraph process (with parameters $\lambda $ and $c$)
independent from the Brownian motions $B_{k},$ $k=2,\ldots,n+1$, we have a
similar result. The distribution $u_{1/2^{n}}$ of (\ref{tel.1}) is a
solution to
\[
\cases{
\dfrac{\partial ^{2/2^{n}}u}{\partial t^{2/2^{n}}}+2\lambda
\dfrac{\partial ^{1/2^{n}}u}{\partial t^{1/2^{n}}}=c^{2}\dfrac{%
\partial ^{2}u}{\partial x^{2}}, \cr
u(x,0)=\delta (x),
}
\qquad x\in \mathbb{R}, t>0
\]
and its characteristic function is equal to
\begin{eqnarray}\label{tel.2}
&&\int_{-\infty }^{+\infty }e^{i\beta x}u_{1/2^{n}}(x,t)\,dx \nonumber \\
&&\qquad=\frac{1}{2}\biggl[ \biggl( 1+\frac{\lambda }{\sqrt{\lambda ^{2}-c^{2}\beta
^{2}}}\biggr) E_{1/2^{n},1}( \eta _{1}t^{1/2^{n}
})\\
&&\hspace*{45pt}
{} +\biggl( 1-\frac{\lambda }{\sqrt{\lambda ^{2}-c^{2}\beta ^{2}}}%
\biggr) E_{1/2^{n},1}( \eta _{2}t^{1/2^{n}}) %
\biggr] ,  \nonumber
\end{eqnarray}%
where $\eta _{1}=-\lambda +\sqrt{\lambda ^{2}-c^{2}\beta ^{2}}$ and $\eta
_{2}=-\lambda -\sqrt{\lambda ^{2}-c^{2}\beta ^{2}}$ [see Orsingher and
Beghin (\citeyear{OrsBeg2004}), formula (2.7), for $\alpha =1/2^{n}$].

For $n\rightarrow \infty $ we get from (\ref{tel.2}) that%
\begin{equation}
\lim_{n\rightarrow \infty }\int_{-\infty }^{+\infty }e^{i\beta x}u_{1/
2^{n}}(x,t)\,dx=\frac{1+2\lambda }{1+2\lambda +c^{2}\beta ^{2}},
\label{tel.4}
\end{equation}%
which is the characteristic function of the bilateral exponential random
variable, with density%
\begin{equation}
f(x)=\frac{\sqrt{1+2\lambda }}{2c}e^{-|x|\sqrt{1+2\lambda }/c
},\qquad x\in \mathbb{R}.  \label{tel.3}
\end{equation}

Clearly, for $\lambda =0$ and $c=1/2$, (\ref{tel.3}) reduces to (\ref{lim.1}%
) and (\ref{tel.4}) coincides with (\ref{lim.4}).
\end{remark}

\section{The explicit solution of the fractional diffusion equation for $%
\nu =1/3$, $\nu =2/3$ and $\nu =4/3$}\label{s4}

In some special cases it is possible to present the solutions of the
fractional equations (\ref{equation1}) in a more attractive fashion. This is
the case for $\nu =\frac{2}{3}$. The explicit form of $u_{2/3}(x,t)$
is given in the next theorem, in terms of Airy functions.

By combining this result with the relationship given in Theorem \ref{thm1}, $%
u_{1/3}(x,t)$ can be represented consequently in an interesting form.

\begin{thm}\label{thm4.1}
The solution to
\begin{equation}\label{tre.3}
\cases{
\dfrac{\partial ^{2/3}u}{\partial t^{2/3}}=\lambda ^{2}\dfrac{\partial ^{2}u}{%
\partial x^{2}}, \cr
u(x,0)=\delta (x),
}
\qquad x\in \mathbb{R}, t>0
\end{equation}%
can be represented as
\begin{equation}\label{tre.4}
u_{2/3}(x,t)=\frac{3}{2}\frac{1}{\lambda \sqrt[3]{3t}}Ai\biggl( \frac{%
|x|}{\lambda \sqrt[3]{3t}}\biggr) ,
\end{equation}%
where
\begin{eqnarray}\label{tre.5}
Ai(w) &=&\frac{1}{\pi }\int_{0}^{+\infty }\cos \biggl( \alpha w+\frac{\alpha
^{3}}{3}\biggr) \,d\alpha  \nonumber\\[-8pt]\\[-8pt]
&=&\frac{w^{1/2}}{3}\biggl[ I_{-1/3}\biggl( \frac{2w^{3/2}}{3}\biggr)
-I_{1/3}\biggl( \frac{2w^{3/2}}{3}\biggr) \biggr]  \nonumber
\end{eqnarray}%
is the Airy function and $I_{\nu }$ denotes the Bessel
function of imaginary argument of order $\nu $.
\end{thm}

\begin{pf}
From (\ref{equation2}) we readily have that
\begin{eqnarray}\label{tre.6}
u_{2/3}(x,t) &=&\frac{1}{2\lambda t^{1/3}}\sum_{k=0}^{\infty }\frac{%
( -|x|/(\lambda t^{1/3})) ^{k}}{k!\Gamma (1-(k+1)/3)}
\nonumber\\[-8pt]\\[-8pt]
&=&\frac{1}{2\pi \lambda t^{1/3}}\sum_{k=0}^{\infty }\frac{( -|x|%
/(\lambda t^{1/3})) ^{k}\Gamma ((k+1)/3)\sin (\pi (k+1)/3)
}{k!}.  \nonumber
\end{eqnarray}

By direct inspection the following identity is proven to hold:%
\begin{equation}
\sin \frac{\pi (k+1)}{3}=(-1)^{k}\sin \frac{2\pi (k+1)}{3}  \label{tre.7}
\end{equation}
and, by inserting this into (\ref{tre.6}), we get that%
\begin{equation}  \label{tre.8}
\qquad u_{2/3}(x,t)=\frac{1}{2\pi \lambda t^{1/3}}\sum_{k=0}^{\infty }\frac{%
( |x|/(\lambda t^{1/3})) ^{k}\Gamma ((k+1)/3)\sin
(2\pi (k+1)/3)}{k!}.
\end{equation}

We note that, from (\ref{tre.5}), for all $|w|<\infty $,%
\begin{eqnarray}\label{tre.9}
Ai(w)
&=&\frac{w^{1/2}}{3}\biggl[ I_{-1/3}\biggl( \frac{2w^{3/2}}{3}\biggr)
-I_{1/3}\biggl( \frac{2w^{3/2}}{3}\biggr) \biggr]  \nonumber \\
&=&\frac{w^{1/2}}{3}\Biggl[ \sum_{k=0}^{\infty }\biggl( \frac{w^{3/2}}{3}%
\biggr) ^{2k-1/3}\frac{1}{k!\Gamma (k-1/3+1)}\nonumber\\
&&\hspace*{28pt}
{}-\sum_{k=0}^{\infty }\biggl( \frac{w^{3/2}}{3}\biggr) ^{2k+1/3}\frac{1%
}{k!\Gamma (k+1/3+1)}\Biggr]   \\
&=&\sum_{k=0}^{\infty }\frac{w^{3k}}{3^{2k+2/3}}\frac{1}{k!\Gamma (k+%
2/3)}-\sum_{k=0}^{\infty }\frac{w^{3k}}{3^{2k+4/3}}\frac{1}{%
k!\Gamma (k+4/3)}  \nonumber \\
&=&\frac{2}{3^{7/6}}\sum_{k=0}^{\infty }\biggl( \frac{w}{3^{2/3}}%
\biggr) ^{k}\frac{\sin (2\pi (k+1)/3)}{\Gamma ((k+2)/3)\Gamma (%
(k+3)/3)}.  \nonumber
\end{eqnarray}

The last step can be justified by taking $k=3m, 3m+1$ and $3m+2$. While for $%
k=3m+2$ the last term in (\ref{tre.9}) is equal to zero, in the other two
cases the two series are obtained.

The triplication formula of the Gamma function [see Lebedev (\citeyear{Leb1972}), page 14],
that is,%
\begin{equation}
\Gamma (z)\Gamma \biggl( z+\frac{1}{3}\biggr) \Gamma\biggl( z+\frac{2}{3}%
\biggr) =\frac{2\pi }{3^{3z-1/2}}\Gamma (3z),  \label{tripl}
\end{equation}%
for $z=\frac{k+1}{3}$ yields%
\begin{equation}
\Gamma \biggl( \frac{k+2}{3}\biggr) \Gamma \biggl( \frac{k+3}{3}\biggr) =\frac{%
2\pi }{3^{k+1/2}}\frac{\Gamma (k+1)}{\Gamma ( (k+1)/3
) }.  \label{tre.10}
\end{equation}

From (\ref{tre.10}) we have that%
\begin{equation}\label{tre.11}
Ai(w)=\frac{3^{-2/3}}{\pi }\sum_{k=0}^{\infty }( 3^{1/3
}w) ^{k}\frac{\sin (2\pi (k+1)/3)}{k!}\Gamma \biggl( \frac{k+1}{3}%
\biggr) ,
\end{equation}%
and (\ref{tre.4}) easily follows by comparing (\ref{tre.11}) and (\ref%
{tre.8}).
\end{pf}

\begin{figure}
\begin{tabular}{@{}l@{}}
\multicolumn{1}{@{}c@{}}{
\includegraphics{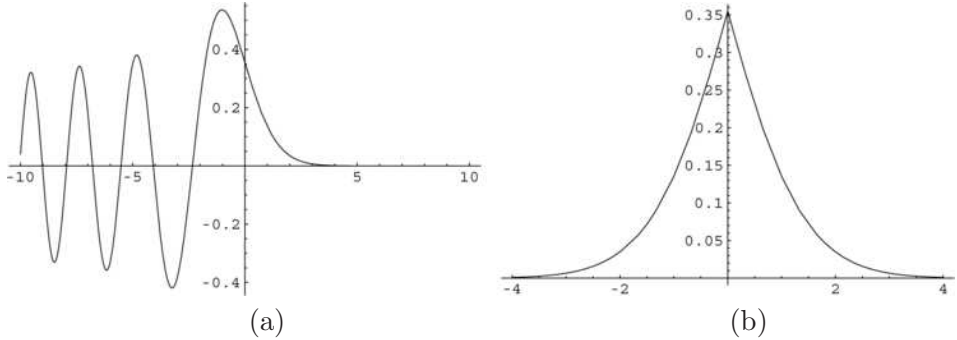}
}\\
\hspace*{32mm}(a)\hspace*{59mm}(b)
\end{tabular}
\caption{The Airy function and the function $Ai(|x|)$.}\label{fig1}
\end{figure}

\begin{remark}\label{rem4.1}
The expression of $u_{2/3}(x,t)$ obtained in the previous theorem
can be recognized (up to the factor $3/2$) as the solution of the
third-order heat-type equation
\begin{equation}\label{tre.13}
\cases{
\displaystyle\frac{\partial v}{\partial t}=-\lambda ^{3}\frac{\partial ^{3}v}{\partial
y^{3}}, \cr
v(y,0)=\delta (y),
}\qquad y\in \mathbb{R}, t>0,
\end{equation}%
evaluated at $y=|x|$. Since $Ai(y)$, for $y>0$, is positive-valued [see
Figure \ref{fig1}(a)] and the function (\ref{tre.4}) integrates to one (as we
show below), $u_{2/3}(x,t)$ is a true probability distribution:%
\begin{eqnarray*}
&&\int_{-\infty }^{+\infty }u_{2/3}(x,t)\,dx \\
&&\qquad=\frac{3}{2}\Biggl[ \int_{0}^{+\infty }\frac{1}{\lambda \sqrt[3]{3t}}%
Ai\biggl( \frac{x}{\lambda \sqrt[3]{3t}}\biggr)\, dx+\int_{-\infty }^{0}\frac{1%
}{\lambda \sqrt[3]{3t}}Ai\biggl( -\frac{x}{\lambda \sqrt[3]{3t}}\biggr) \,dx%
\Biggr] \\
&&\qquad=2\frac{3}{2}\int_{0}^{+\infty }Ai( y) \,dy=1,
\end{eqnarray*}%
where the last step follows by noting that $\int_{0}^{+\infty }Ai(
y) \,dy=1/3$; see Nikitin and Orsingher (\citeyear{NikOrs2000}).

Therefore we can think of $u_{2/3}(x,t)$ as the probability law of a
process $A(t),t>0,$ whose distribution at time $t$ is obtained from the
solution $v(x,t)$ of equation (\ref{tre.13}), as follows:%
\[
u_{2/3}(x,t)=\tfrac{3}{2}v(|x|,t).
\]
\end{remark}

\begin{remark}\label{rem4.2}
For the case $\nu =\frac{1}{3}$ the solution $u_{1/3}(x,t)$ to (\ref%
{equation1}) can be written, thanks to the relationship (\ref{vologda}), as%
\begin{eqnarray}\label{tre.14}
u_{1/3}(x,t) &=&\frac{1}{\sqrt{\pi t}}\int_{0}^{\infty }e^{-
z^{2}/(4t)}u_{2/3}(x,z)\,dz \\ 
&=&\frac{1}{\sqrt{\pi t}}\int_{0}^{\infty }e^{-z^{2}/(4t)}\frac{3^{2/3}%
}{2\lambda z^{1/3}}Ai\biggl( \frac{|x|}{\lambda \sqrt[3]{3z}}\biggr) \,dz.
\nonumber
\end{eqnarray}

We can represent (\ref{tre.14}) as the distribution of the process%
\[
J_{1/3}(t)=A(|B(t)|),\qquad t>0,
\]
with $A$ and $B$ independent. The results (\ref{tre.4}) and (\ref{tre.14})
show that the solutions $u_{2/3}(x,t)$ and $u_{1/3}(x,t)$
are both unimodal with maximum at $x=0$; see Figure~\ref{fig1}(b). This is in
accordance with the general result that, for $0<\nu \leq 1$, the solutions
to the fractional equation (\ref{equation1}) have a unique maximal point at $%
x=0$.
\end{remark}

We consider now the case $\nu =4/3$, which is qualitatively different from
those dealt with so far, because the solutions of fractional equations of
order $1<\nu <2$ display a substantially different behavior.

\begin{thm}\label{thm4.2}
The solution to
\begin{equation} \label{qua.3}
\cases{
\displaystyle\frac{\partial ^{4/3}u}{\partial t^{4/3}}=\lambda ^{2}\frac{\partial ^{2}u}{%
\partial x^{2}}, \cr
\displaystyle u(x,0)=\delta (x), \cr
\displaystyle u_{t}(x,0)=0,
}
\qquad x\in \mathbb{R}, t>0,
\end{equation}%
is given by
\begin{equation}\label{qua3.bis}
\qquad\quad u_{4/3}(x,t)=\frac{1}{\lambda \sqrt{\pi }}\biggl( \frac{3}{4t}\biggr)
^{2/3}\int_{0}^{+\infty }e^{-w}w^{-1/6}Ai\biggl( -\frac{|x|}{%
\lambda }\biggl( \frac{2}{t}\sqrt{\frac{w}{3}}\biggr) ^{2/3}\biggr)\,dw.
\end{equation}
\end{thm}

\begin{pf}
From (\ref{equation2}) we have that
\begin{eqnarray} \label{qua.4}
u_{4/3}(x,t) &=&\frac{1}{2\lambda t^{2/3}}\sum_{k=0}^{\infty }\biggl(
-\frac{|x|}{\lambda t^{2/3}}\biggr) ^{k}\frac{1}{k!\Gamma (1-2/3
(k+1))}  \nonumber\\[-8pt]\\[-8pt]
&=&\frac{1}{2\lambda \pi t^{2/3}}\sum_{k=0}^{\infty }\biggl( -\frac{|x|}{%
\lambda t^{2/3}}\biggr) ^{k}\frac{\Gamma (2/3(k+1))\sin (2\pi
(k+1)/3)}{k!}.  \nonumber
\end{eqnarray}

By means of the duplication formula for the Gamma function we have that%
\[
\Gamma \biggl( \frac{1}{3}(k+1)+\frac{1}{2}\biggr) =\frac{\sqrt{\pi }2^{1-
2/3(k+1)}\Gamma (2/3(k+1))}{\Gamma ((k+1)/3)},
\]
and therefore $u_{4/3}(x,t)$ can be rewritten as%
\begin{eqnarray}\label{qua.5}
u_{4/3}(x,t)
&=&\frac{1}{2\lambda \pi \sqrt{\pi }2^{1/3}t^{2/3}}\nonumber\\
&&{}\times\sum_{k=0}^{\infty
}\biggl( -\frac{2^{2/3}|x|}{\lambda t^{2/3}}\biggr) ^{k}\frac{\sin (2\pi
(k+1)/3)}{k!}\Gamma \biggl( \frac{k+1}{3}\biggr) \Gamma \biggl( \frac{k+1}{3}+%
\frac{1}{2}\biggr)   \\
&=&\frac{1}{2\lambda \pi \sqrt{\pi }2^{1/3}t^{2/3}}\nonumber\\
&&{}\times\sum_{k=0}^{\infty
}\int_{0}^{+\infty }e^{-w}w^{1/3(k+1)+1/2-1}\nonumber \\
&&\phantom{{}\times\sum_{k=0}^{\infty
}\int_{0}^{+\infty }}
{}\times\Gamma \biggl(
\frac{k+1}{3}\biggr) \biggl( -\frac{|x|}{\lambda }\biggl( \frac{2}{t}\biggr) ^{%
2/3}\biggr) ^{k}\frac{\sin( 2\pi (k+1)/3)}{k!}\,dw  \nonumber \\
&=&\frac{1}{2\lambda \pi \sqrt{\pi }2^{1/3}t^{2/3}}\nonumber\\
&&{}\times\sum_{k=0}^{\infty
}\int_{0}^{+\infty }e^{-w}w^{1/2-2/3}\biggl( -\frac{|x|}{%
\lambda }\biggl( \frac{2}{t}\biggr) ^{2/3}w^{1/3}\biggr)
^{k}\nonumber\\
&&\phantom{{}\times\sum_{k=0}^{\infty
}\int_{0}^{+\infty }}
{}\times
\frac{\sin (2\pi (k+1)/3)}{k!}\Gamma \biggl( \frac{k+1}{3}\biggr) \,dw \nonumber\\
&=&[ \mbox{by (\ref{tre.11})}] \nonumber\\
&=&\frac{3^{2/3}}{2\lambda \sqrt{\pi }2^{1/3}t^{2/3}}\int_{0}^{+\infty
}e^{-w}w^{-1/6}Ai\biggl( -\frac{|x|}{\lambda }\biggl( \frac{2}{t}\sqrt{%
\frac{w}{3}}\biggr) ^{2/3}\biggr) \,dw.\nonumber
\end{eqnarray}
\upqed
\end{pf}

We can show that $\int_{-\infty }^{+\infty }u_{4/3}(x,t)\,dx=1$.
Indeed, from (\ref{qua3.bis}) we have that%
\begin{eqnarray*}
&&\frac{1}{\lambda \sqrt{\pi }}\biggl( \frac{3}{4t}\biggr) ^{2/3}
\int_{0}^{+\infty }e^{-w}w^{-1/6}\int_{-\infty }^{+\infty }Ai\biggl(
-\frac{|x|}{\lambda }\biggl( \frac{2}{t}\sqrt{\frac{w}{3}}\biggr) ^{2/3}\biggr) \,dx\,dw \\
&&\qquad=\frac{2}{\lambda \sqrt{\pi }}\biggl( \frac{3}{4t}\biggr) ^{2/3}
\int_{0}^{+\infty }e^{-w}w^{-1/6}\int_{0}^{+\infty }Ai\biggl( -\frac{%
x}{\lambda }\biggl( \frac{2}{t}\sqrt{\frac{w}{3}}\biggr) ^{2/3}\biggr) \,dx\,dw\\
&&\qquad=\biggl[ \mbox{by the substitution }y=-\frac{x}{\lambda }\biggl( \frac{2}{t}%
\sqrt{\frac{w}{3}}\biggr) ^{2/3}\biggr] \\
&&\qquad=\frac{2}{\lambda \sqrt{\pi }}\biggl( \frac{3}{4t}\biggr) ^{2/3}\lambda
\biggl( \frac{2}{t}\sqrt{\frac{1}{3}}\biggr) ^{-2/3}\int_{0}^{+\infty }e^{-w}w^{-1/2}\int_{-\infty }^{0}Ai(
y) \,dy\,dw \\
&&\qquad=\frac{1}{\sqrt{\pi }}\frac{2}{3}2\biggl( \frac{3}{4}\biggr) ^{2/3}%
\frac{2^{-2/3}}{3^{-1/3}}\int_{0}^{+\infty }e^{-w}w^{-1/2}\,dw=1.
\end{eqnarray*}

\begin{remark}\label{rem4.3}
In view of Theorem \ref{thm1} we have the following representation for
$u_{2/3}(x,t)$, which is alternative to (\ref{tre.4}):%
\begin{eqnarray}\label{qua.6}
u_{2/3}(x,t) &=&\frac{1}{\sqrt{\pi t}}\int_{0}^{\infty }e^{-
z^{2}/(4t)}u_{4/3}(x,z)\,dz  \\
&=&\frac{3^{2/3}}{2\lambda \sqrt{\pi }2^{1/3}}\int_{0}^{+\infty }\frac{e^{-%
z^{2}/(4t)}}{z^{2/3}\sqrt{\pi t}}\,dz\nonumber\\
&&{}\times\int_{0}^{+\infty }e^{-w}w^{-1/6}Ai\biggl( -\frac{|x|}{\lambda }\biggl( \frac{2}{z}\sqrt{\frac{w}{3}}%
\biggr) ^{2/3}\biggr) \,dw.  \nonumber
\end{eqnarray}

By inserting (\ref{tre.4}) into the left-hand side of (\ref{qua.6}) we
obtain that%
\begin{eqnarray}\label{qua.ai}
&&\frac{3}{2\lambda \sqrt[3]{3t}}Ai\biggl( \frac{|x|}{\lambda \sqrt[3]{3t}}\biggr)\nonumber   \\
&&\qquad=\frac{3^{2/3}}{2\lambda \sqrt{\pi }2^{1/3}}\int_{0}^{+\infty }\frac{e^{-%
z^{2}/(4t)}}{\sqrt{\pi t}z^{2/3}}\,dz\nonumber\\
&&\quad\qquad{}\times\int_{0}^{+\infty }e^{-w}w^{-1/6}Ai\biggl( -\frac{|x|}{\lambda }\biggl( \frac{2}{z}\sqrt{\frac{w}{3}}%
\biggr) ^{2/3}\biggr)\, dw  \nonumber \\
&&\qquad=\bigl[ \mbox{by the substitution }s=\sqrt[3]{2^{2}twz^{-2}}\bigr]  \nonumber
\\
&&\qquad=\frac{3^{2/3}}{2\lambda \sqrt{\pi }2^{1/3}}\int_{0}^{+\infty }\frac{e^{-%
z^{2}/(4t)}}{\sqrt{\pi t}z^{2/3}}\frac{3z^{5/3}}{2^{5/3
}t^{5/6}}s^{3/2}\,dz\\
&&\quad\qquad{}\times\int_{0}^{+\infty }e^{-z^{2}s^{3}/(4t)}Ai\biggl( -\frac{|x|s}{\lambda \sqrt[3]{3t}}\biggr) \,ds  \nonumber \\
&&\qquad=\frac{3^{5/3}}{2^{3}\lambda \pi t^{4/3}}\int_{0}^{+\infty }Ai\biggl( -%
\frac{|x|s}{\lambda \sqrt[3]{3t}}\biggr) s^{3/2}\,ds\int_{0}^{+\infty
}ze^{-z^{2}(1+s^{3})/(4t)}\,dz  \nonumber \\
&&\qquad=\frac{3^{5/3}}{2^{2}\lambda \pi t^{1/3}}\int_{0}^{+\infty }\frac{s^{
3/2}}{1+s^{3}}Ai\biggl( -\frac{|x|s}{\lambda \sqrt[3]{3t}}\biggr) \,ds  \nonumber
\\
&&\qquad=\frac{3^{2/3}}{2\lambda t^{1/3}}\int_{0}^{+\infty }\Pr \{ |B(
T_{0}) |\in ds\} Ai\biggl( -\frac{|x|s}{\lambda \sqrt[3]{3t}}%
\biggr) ,  \nonumber
\end{eqnarray}%
where
\[
\Pr \{ |B( T_{0}) |\in ds\} =\frac{3}{2\pi }\frac{s^{%
3/2}}{1+s^{3}}\,ds,\qquad s>0,
\]
is the McKean law representing the distribution of the position of a
Brownian motion $B$ at the instant
\[
T_{0}=\inf \biggl\{ t>0\dvtx1+\int_{0}^{t}B(s)\,ds=0\biggr\};
\]
see McKean (\citeyear{McK1963}).

By setting $y=\frac{|x|}{\lambda \sqrt[3]{3t}}$ in (\ref{qua.ai}) and
performing some simplifications we get%
\begin{equation}
Ai(|y|)=\int_{0}^{+\infty }\Pr \{ |B( T_{0}) |\in ds\}
Ai( -|y|s) ,\qquad y\in \mathbb{R}.  \label{airy6}
\end{equation}

Formula (\ref{airy6}) shows an interesting property of Airy functions: The
value of the exponentially decreasing part of $Ai(|y|)$ can be obtained by
averaging its oscillating component $Ai( -|y|s) $\ with the
well-known density of $|B( T_{0}) |$ (see Figure~\ref{fig1}).
\end{remark}

\begin{remark}\label{rem4.4}
The solution $u_{4/3}(x,t)$ can also be expressed in terms of a
stable density of order $\frac{3}{2}.$ Indeed, by using the representation
of the stable density below%
\begin{equation}
\qquad\quad p_{\alpha }(x;\gamma ,\eta )=\frac{1}{2\pi }\int_{-\infty }^{+\infty
}e^{-i\beta x}\exp \bigl\{ -\eta |\beta |^{\alpha }e^{-i\pi \gamma /2
\beta /|\beta |}\bigr\}\, d\beta ,\qquad \alpha \neq 1,  \label{stab}
\end{equation}%
we know that for $\alpha \in ( 1,2) $, $\eta =1$ and for $x>0$
the following series representation holds true:%
\begin{equation}
p_{\alpha }(x;\gamma ,1)=\frac{1}{\pi }\sum_{k=1}^{\infty }( -x)
^{k-1}\frac{\sin (k\pi (\gamma +\alpha )/(2\alpha))}{k!}\Gamma \biggl( 1+%
\frac{k}{\alpha }\biggr);  \label{qua.7}
\end{equation}%
see formula (6.9), page 583 of Feller (\citeyear{Fel1971}) (up to some corrections) and
Lukacs (\citeyear{Luk1969}).

For $\alpha =\frac{3}{2}$ and $\gamma =\frac{1}{2}$ formula (\ref{qua.7})
reads%
\begin{eqnarray} \label{qua.8}
p_{3/2}\biggl(x;\frac{1}{2},1\biggr) &=&\frac{1}{\pi }\sum_{r=0}^{\infty }(
-x) ^{r}\frac{\sin \{ (r+1)2/3\pi \} }{(r+1)!}%
\Gamma \biggl( 1+\frac{2}{3}(r+1)\biggr)  \nonumber\\[-8pt]\\[-8pt]
&=&\frac{2}{3}\frac{1}{\pi }\sum_{r=0}^{\infty }( -x) ^{r}\frac{%
\sin \{ (r+1)2/3\pi \} }{r!}\Gamma \biggl( \frac{2}{3}%
(r+1)\biggr) .  \nonumber
\end{eqnarray}

If we compare (\ref{qua.4}) with (\ref{qua.8}) we get that%
\begin{equation}
u_{4/3}(x,t)=\frac{3}{2}\frac{1}{2\lambda t^{2/3}}p_{3/2
}\biggl( \frac{|x|}{\lambda t^{2/3}};\frac{1}{2},1\biggr) .  \label{qua.9}
\end{equation}

A different proof of the relationship between stable laws and the solutions
of fractional diffusion equations, based on the inversion of the Fourier
transform, can be found in Fujita (\citeyear{Fuj1990}).

Formula (\ref{qua.9}) proves the nonnegativity of the expression (\ref%
{qua.5}), as a function of~$x$.
\end{remark}

\section{Some generalizations of the previous results}\label{s5}

In this section we present some generalizations of the results of Sections \ref{s2}
and \ref{s4}.

We start by giving a relationship between the solutions $u_{\nu }$ and $%
u_{m\nu }$, $m\geq 3$, and obtain some explicit expressions for $m=3$. In
this case the interpretation of $u_{2/3^{n}}$ as the distribution of
compositions of different types of processes is possible. Also in this case
we encounter processes with a random time which possesses a branching
structure (depending on $n$).

We now state a general result which is alternative to (\ref{vologda}) and
permits us to exploit the explicit expression of $u_{\nu }(x,t)$.

\begin{thm}\label{thm5.1}
The solution to the initial value
problem (\ref{equation1})--(\ref{condition 1}), for $0<\nu \leq 2/3$, can be represented as
\begin{equation}
u_{\nu }(x,t)=\frac{3}{2\pi \sqrt{t}}\int_{0}^{+\infty }\int_{0}^{+\infty
}se^{-(s^{3}+v^{3})/(3\sqrt{3t})}u_{3\nu }( x,sv) \,ds\,dv,
\label{cin.2}
\end{equation}
where $u_{3\nu }(x,z)$ is the solution to
\begin{equation}
\cases{
\displaystyle\frac{\partial ^{3\nu }u}{\partial z^{3\nu }}=\lambda ^{2}\frac{\partial
^{2}u}{\partial x^{2}}, \cr
u(x,0)=\delta (x),
}
\qquad x\in \mathbb{R},z>0,\ 0<\nu \leq \frac{1}{3},
\label{cin.3}
\end{equation}%
and%
\begin{equation}
\cases{
\displaystyle\frac{\partial ^{3\nu }u}{\partial z^{3\nu }}=\lambda ^{2}\frac{\partial
^{2}u}{\partial x^{2}}, \cr
u(x,0)=\delta (x), \cr
u_{t}(x,0)=0,
}
\qquad x\in \mathbb{R}, z>0,\  \frac{1}{3}<\nu <\frac{2}{3}.
\label{cin.bis}
\end{equation}
\end{thm}

\begin{pf}
In view of the triplication formula (\ref{tripl}),
for $z=\frac{1}{3}-\frac{\nu (k+1)}{2}$, we have that
\begin{eqnarray}\label{cin.4}
u_{\nu }(x,t)
&=&\frac{1}{2\lambda t^{\nu /2}}\sum_{k=0}^{\infty }\frac{( -|x|/
(\lambda t^{\nu /2})) ^{k}}{k!\Gamma (1-\nu (k+1)/2)}  \nonumber \\
&=&\frac{1}{2\lambda 2\pi t^{\nu /2}}\nonumber\\
&&{}\times\sum_{k=0}^{\infty }\biggl( -\frac{|x|}{%
\lambda t^{\nu /2}}\biggr) ^{k}\nonumber\\
&&\phantom{{}\times\sum_{k=0}^{\infty }}{}\times\frac{3^{1-3/2\nu (k+1)-1/2
}\Gamma (2/3-\nu (k+1)/2)\Gamma (1/3-\nu (k+1)/2)}{k!\Gamma (1-3\nu (k+1)/2)}  \nonumber \\
&=&\frac{\sqrt{3}}{2^{2}3^{3/2\nu }\lambda \pi t^{\nu /2}}\nonumber\\
&&{}\times\int_{0}^{+\infty }\int_{0}^{+\infty }e^{-w-z}w^{-\nu /2-1/3
}z^{-\nu /2-2/3}\nonumber\\
&&\phantom{{}\times\int_{0}^{+\infty }\int_{0}^{+\infty }}
{}\times\sum_{k=0}^{\infty }\biggl( -\frac{|x|}{%
\lambda (\sqrt[3]{3^{3}wzt})^{3\nu /2}}\biggr) ^{k}\frac{dw\,dz}{k!\Gamma (1-%
3\nu (k+1)/2)}  \nonumber \\
&=&\frac{\sqrt{3}}{2\pi 3^{3/2\nu }t^{\nu /2}}\int_{0}^{+\infty
}\int_{0}^{+\infty }e^{-w-z}w^{-\nu /2-1/3}z^{-\nu /2
-2/3}\nonumber\\
&&\phantom{\frac{\sqrt{3}}{2\pi 3^{3/2\nu }t^{\nu /2}}\int_{0}^{+\infty
}\int_{0}^{+\infty }}
{}\times\bigl(\sqrt[3]{3^{3}wzt}\bigr)^{3\nu /2}u_{3\nu }\bigl(x,\sqrt[3]{%
3^{3}wzt}\bigr)\,dw\,dz   \\
&=&\frac{\sqrt{3}(3^{3}t)^{\nu /2}}{2\pi 3^{3/2\nu }t^{\nu
/2}}\int_{0}^{+\infty }\int_{0}^{+\infty }e^{-w-z}w^{-\nu /2-1%
/3}z^{-\nu /2-2/3}\nonumber\\
&&\phantom{\frac{\sqrt{3}(3^{3}t)^{\nu /2}}{2\pi 3^{3/2\nu }t^{\nu
/2}}\int_{0}^{+\infty }\int_{0}^{+\infty }}
{}\times(wz)^{\nu /2}u_{3\nu }\bigl(x,\sqrt[3]%
{3^{3}wzt}\bigr)\,dw\,dz \nonumber  \\
&=&\frac{\sqrt{3}}{2\pi }\int_{0}^{+\infty }\int_{0}^{+\infty }e^{-w-z}w^{-%
1/3}z^{-2/3}u_{3\nu }\bigl(x,3\sqrt[3]{wzt}\bigr)\,dw\,dz,  \nonumber
\end{eqnarray}%
which reduces to (\ref{cin.2}), after the change of variables%
\[
\cases{
s=\sqrt{3}\sqrt[3]{w}\sqrt[3]{t^{1/2}}, \cr
v=\sqrt{3}\sqrt[3]{z}\sqrt[3]{t^{1/2}}.
}
\]
\upqed
\end{pf}

It can be easily checked that, also in this form, the solution integrates to
one. By using the last expression in (\ref{cin.4}) we get%
\begin{eqnarray*}
&&\int_{-\infty }^{+\infty }u_{\nu }(x,t)\,dx\\
&&\qquad=\frac{\sqrt{3}}{2\pi }\int_{0}^{+\infty }\int_{0}^{+\infty }e^{-w-z}w^{-%
1/3}z^{-2/3}\int_{-\infty }^{+\infty }u_{3\nu }\bigl(x,3\sqrt[3]{%
wzt}\bigr)\,dx\,dw\,dz \\
&&\qquad=\frac{\sqrt{3}}{2\pi }\int_{0}^{+\infty }e^{-w}w^{-1+2/3
}\,dw\int_{0}^{+\infty }e^{-z}z^{-1+1/3}\,dz \\
&&\qquad=\frac{\sqrt{3}}{2\pi }\Gamma \biggl( \frac{2}{3}\biggr) \Gamma \biggl(
\frac{1}{3}\biggr) =1,
\end{eqnarray*}%
since, by the triplication formula for $z=1/3$, it is $\Gamma ( \frac{2%
}{3}) \Gamma ( \frac{1}{3}) =2\pi /\sqrt{3}.$

\begin{remark}\label{rem5.1}
By using the previous result it is possible to obtain alternative forms for
the solution to the initial value problem for $\nu =1/3$ and for $\nu =2/9.$
Indeed, in the first case it is%
\begin{eqnarray}\label{cin.5}
u_{1/3}(x,t) &=&\frac{3}{2\pi \sqrt{t}}\int_{0}^{+\infty }\int_{0}^{+\infty
}se^{-(s^{3}+v^{3})/(3\sqrt{3t})}u_{1}( x,sv) \,ds\,dv
\nonumber\\[-8pt]\\[-8pt]
&=&\frac{3}{2\pi \sqrt{t}}\int_{0}^{+\infty }\int_{0}^{+\infty }se^{-
(s^{3}+v^{3})/(3\sqrt{3t})}\frac{e^{-x^{2}/(4\lambda ^{2}(sv))}}{%
2\lambda \sqrt{\pi sv}}\,ds\,dv.  \nonumber
\end{eqnarray}

The relationship (\ref{cin.5}) shows that $u_{1/3}$ can be interpreted as
the distribution of a Brownian motion (with infinitesimal variance $2\lambda
^{2}$) at a random time $G_{1}(t)\cdot G_{2}(t)$, that is,%
\begin{equation}
J_{1/3}(t)=B[ G_{1}(t)\cdot G_{2}(t)] ,  \label{cin.5bis}
\end{equation}%
where $(G_{1}(t),G_{2}(t))$ possesses joint density%
\begin{equation}
p_{(G_{1}(t),G_{2}(t))}(s,v)=\frac{3}{2\pi \sqrt{t}}se^{-(s^{3}+v^{3})
/(3\sqrt{3t})},\qquad s>0,v>0.  \label{cin.6}
\end{equation}

This result corresponds to (\ref{iterated Brownian motion}), for $\nu =1/3$
and it represents a counterpart of result (\ref{iterated}) with the
reflecting Brownian motion replaced by the product $G_{1}(t)\cdot G_{2}(t),$
with joint distribution given in (\ref{cin.6}).

In the case $\nu =2/3^{2}$, from (\ref{cin.2}) we have that%
\begin{equation}
u_{2/3^{2}}(x,t)=\frac{3}{2\pi \sqrt{t}}\int_{0}^{+\infty
}\int_{0}^{+\infty }se^{-(s^{3}+v^{3})/(3\sqrt{3t})}u_{2/3}( x,sv) \,ds\,dv  \label{nuova}
\end{equation}%
and this suggests that we interpret $u_{2/3^{2}}(x,t)$ as the
distribution of the process%
\begin{equation}
J_{2/3^{2}}(t)=A[ G_{1}(t)\cdot G_{2}(t)] .  \label{airy}
\end{equation}

The process (\ref{airy}) is analogous to (\ref{cin.5bis}) with the role of
Brownian motion played by the process $A$.

Analogously to (\ref{airy}), for $\nu =2/3^{3}$, we get%
\[
J_{2/3^{3}}(t)=A\{ G_{1}[ G_{1}(t)\cdot G_{2}(t)]
\cdot G_{2}[ G_{1}(t)\cdot G_{2}(t)] \}
\]
which has distribution coinciding with%
\begin{eqnarray*}
&&u_{2/3^{3}}(x,t)\\
&&\qquad=\frac{3}{2\pi \sqrt{t}}\int_{0}^{+\infty }\!\!\int_{0}^{+\infty }se^{-
(s^{3}+v^{3})/(3\sqrt{3t})}\\
&&\quad\qquad{}\times\biggl( \frac{3}{2\pi \sqrt{sv}}\int_{0}^{+\infty
}\!\!\int_{0}^{+\infty }we^{-(w^{3}+z^{3})/(3\sqrt{3sv})}u_{2/3
}( x,zw) \,dz\,dw\biggr) \,ds\,dv,
\end{eqnarray*}%
as an application of (\ref{cin.2}) and (\ref{nuova}) shows.
\end{remark}

The results of Theorem \ref{thm1} and \ref{thm5.1} can be furthermore generalized in order
to relate the solutions $u_{\nu }(x,t)$ with $u_{m\nu }(x,t)$.

\begin{thm}\label{thm5.2}
The solution to equation (\ref{eq4}),
for $\nu \leq 2/m$, $m\geq 1$ can be represented as
\begin{eqnarray}\label{sei.1}
u_{\nu }(x,t)
&=&\frac{m^{(m-1)/2}}{(2\pi )^{(m-1)/2}\sqrt{t}}\nonumber\\
&&
{}\times
\int_{0}^{+\infty }\cdots\int_{0}^{+\infty }e^{-(w_{1}^{m}+\cdots+w_{m-1}^{m})
/\sqrt[m-1]{m^{m}t}}w_{2}\cdots w_{m-1}^{m-2}\\
&&\hspace*{77pt}
{}\times u_{m\nu }(
x,w_{1}w_{2}\cdots w_{m-1}) \,dw_{1}\cdots dw_{m-1}.  \nonumber
\end{eqnarray}
\end{thm}

\begin{pf}
From (\ref{equation2}), by using the multiplication
formula of the Gamma function [see Magnus and Oberhettinger (\citeyear{MagObe1948})], that is,%
\begin{eqnarray*}
&&\Gamma (z)\Gamma \biggl( z+\frac{1}{m}\biggr) \Gamma \biggl( z+\frac{2}{m}%
\biggr) \cdots\Gamma \biggl( z+\frac{m-1}{m}\biggr) \\
&&\qquad=( 2\pi ) ^{(m-1)/2}m^{1/2-mz}\Gamma (mz),
\end{eqnarray*}%
for $z=\frac{1}{m}-\frac{\nu (k+1)}{2}$, we get that%
\begin{eqnarray*}
&&u_{\nu }(x,t)\\
&&\qquad=\frac{1}{2\lambda t^{\nu /2}}\sum_{k=0}^{\infty }\biggl( -\frac{|x|}{%
\lambda t^{\nu /2}}\biggr) ^{k}\frac{1}{k!\Gamma (1-\frac{\nu (k+1)}{2})} \\
&&\qquad=\frac{\sqrt{m}}{2\lambda t^{\nu /2}(2\pi )^{(m-1)/2}}\\
&&\quad\qquad{}\times
\sum_{k=0}^{\infty }\biggl( -\frac{|x|}{\lambda t^{\nu /2}}\biggr) ^{k}
\\
&&\quad\qquad\phantom{{}\times\sum_{k=0}^{\infty }}
{}\times \Gamma \biggl(\frac{1}{m}-\frac{\nu (k+1)}{2}\biggr)\Gamma \biggl( \frac{2}{m}-
\frac{\nu (k+1)}{2}\biggr) \cdots\\
&&\quad\qquad\phantom{{}\times\sum_{k=0}^{\infty }}
{}\times \Gamma \biggl( \frac{m-1}{m}-\frac{\nu (k+1)}{2}\biggr) m^{-m/2\nu (k+1)}\biggl[k!\Gamma \biggl(1-\frac{m}{2}\nu (k+1)\biggr)\biggr]^{-1} \\
&&\qquad=\frac{\sqrt{m}}{2\lambda t^{\nu /2}(2\pi )^{(m-1)/2}m^{m/2\nu }}\\
&&\quad\qquad
{}\times\int_{0}^{+\infty }\cdots\int_{0}^{+\infty }e^{-w_{1}-\cdots-w_{m-1}}w_{1}^{-%
1/m}w_{2}^{-2/m}\cdots w_{m-1}^{-(m-1)/m} \\
&&\quad\qquad\phantom{{}\times\int_{0}^{+\infty }\cdots\int_{0}^{+\infty }}{}\times u_{m\nu }\bigl( x,m\sqrt[m]{w_{1}\cdots w_{m-1}t}\bigr)\,
dw_{1}\cdots dw_{m-1}.
\end{eqnarray*}

By means of the transformation%
\[
z_{j}=\sqrt[m-1]{m}\sqrt[m]{w_{j}}\sqrt[m]{t^{1/(m-1)}},
\]
we finally get (\ref{sei.1}).
\end{pf}

We prove now a general result, valid for any $0<\nu <2$, which gives another
representation for the solution $u_{\nu }=u_{\nu }(x,t)$, alternative to
those presented in the previous sections.

\begin{thm}\label{thm5.3}
The solution to (\ref{equation1}) with
initial condition (\ref{condition 1}) or (\ref{cond1}) has the following form:
\begin{eqnarray}\label{sin}
u_{\nu }(x,t)
&=&\frac{1}{2\pi \lambda t^{\nu /2}}\int_{0}^{+\infty }e^{-w}w^{\nu /2-1}
e^{-|x|w^{\nu /2}/(\lambda t^{\nu /2})\cos ( \nu \pi /2
) }\nonumber\\
&&\phantom{\frac{1}{2\pi \lambda t^{\nu /2}}\int_{0}^{+\infty }}{}\times\sin \biggl( \frac{\nu \pi }{2}-\frac{|x|w^{\nu /2}}{\lambda t^{\nu
/2}}\sin \biggl(\frac{\nu \pi }{2}\biggr)\biggr) \,dw   \\
&=&\frac{1}{\nu \pi }\int_{0}^{+\infty }e^{-|x|y\cos ( \nu \pi /2) -(\lambda y)^{2/\nu }t}\sin \biggl( \frac{\nu \pi }{2}-|x|y\sin \biggl(%
\frac{\nu \pi }{2}\biggr)\biggr) \,dy,  \nonumber
\end{eqnarray}%
for $0<\nu <2$.
\end{thm}

\begin{pf}
By applying the reflection property of the Gamma
function we rewrite (\ref{equation2}) as%
\begin{eqnarray}\label{sin.2}
u_{\nu }(x,t) &=&\frac{1}{2\lambda t^{\nu /2}}\sum_{k=0}^{\infty }\frac{%
( -|x|/(\lambda t^{\nu /2})) ^{k}}{k!\Gamma (1-\nu /2(k+1))} \nonumber\\
&=&\frac{1}{2\pi \lambda t^{\nu /2}}\sum_{k=0}^{\infty }\biggl( -\frac{|x|}{%
\lambda t^{\nu /2}}\biggr) ^{k}\frac{\sin ( \nu \pi /2
(k+1)) }{k!}\Gamma \biggl(\frac{\nu }{2}(k+1)\biggr)  \nonumber \\
&=&\frac{1}{2\pi \lambda t^{\nu /2}}\int_{0}^{+\infty
}e^{-w}\sum_{k=0}^{\infty }\frac{w^{\nu (k+1)/2-1}}{k!}\biggl( -\frac{%
|x|}{\lambda t^{\nu /2}}\biggr) ^{k}\sin \biggl( \frac{\nu \pi }{2}%
(k+1)\biggr) \, dw  \nonumber \\
&=&\frac{1}{2\pi \lambda t^{\nu /2}}\int_{0}^{+\infty }e^{-w}w^{\nu
/2-1}\nonumber\\[-8pt]\\[-8pt]
&&\phantom{\frac{1}{2\pi \lambda t^{\nu /2}}\int_{0}^{+\infty }}
{}\times\sum_{k=0}^{\infty }\biggl[ \biggl( -\frac{|x|w^{\nu /2}}{\lambda
t^{\nu /2}}\biggr) ^{k}\frac{1}{k!}\frac{e^{i\nu \pi (k+1)/2}-e^{-%
i\nu \pi (k+1)/2}}{2i}\biggr] \,dw  \nonumber \\
&=&\frac{1}{2\pi \lambda t^{\nu /2}}\nonumber\\
&&{}\times\int_{0}^{+\infty }e^{-w}\frac{w^{
\nu /2-1}}{2i}\nonumber\\
&&\phantom{{}\times\int_{0}^{+\infty }}
{}\times\bigl[ e^{-|x|w^{\nu /2}e^{i\nu /2\pi} /
(\lambda t^{\nu /2})}e^{i\nu /2\pi }\nonumber\\
&&\hspace*{54pt}
{}-e^{-|x|w^{\nu /2
}e^{-i\nu /2\pi }/(\lambda t^{\nu /2})}e^{-i\nu /2\pi }\bigr] \,dw,  \nonumber
\end{eqnarray}%
which coincides with the first form of (\ref{sin}). The second line can be
obtained by the change of variable $w=(\lambda y)^{2/\nu }t$.
\end{pf}

\begin{remark}\label{rem5.2}
We can check that, for $\nu =1$ (i.e., for the heat equation), the first
expression in (\ref{sin}) reduces to the Gaussian density:
\begin{eqnarray} \label{gaus}
u_{1}(x,t)
&=&\frac{1}{2\pi \lambda t^{1/2}}\int_{0}^{+\infty }e^{-w}w^{1/2
-1}\sin \biggl( \frac{\pi }{2}-\frac{|x|w^{1/2}}{\lambda t^{1/2}}\biggr)\, dw
\nonumber \\
&=&\frac{1}{2\pi \lambda t^{1/2}}\int_{0}^{+\infty }e^{-w}w^{1/2
-1}\cos \biggl( \frac{|x|w^{1/2}}{\lambda t^{1/2}}\biggr) \,dw
\nonumber\\[-8pt]\\[-8pt]
&=&[ w=y^{2}]  \nonumber \\
&=&\frac{1}{\sqrt{4\pi t\lambda ^{2}}}e^{-x^{2}/(4t\lambda ^{2})}.
\nonumber
\end{eqnarray}%
In the last step we used formula 3.896.4, page 514, of Gradshteyn and
Ryzhik (\citeyear{GraRyz1994}). The same check can be done for the second expression in (\ref{sin}).
\end{remark}

An alternative form of (\ref{sin}) can be obtained by means of a double
integration by parts, as follows:

\begin{corollary}\label{coro5.1}
The solution to (\ref{equation1}) with initial condition
(\ref{condition 1}) or (\ref{cond1}) can be rewritten as
\begin{eqnarray} \label{sin.3}
u_{\nu }(x,t)
&=&\frac{1}{\pi \nu |x|}\int_{0}^{+\infty }e^{-w}e^{-|x|w^{\nu /2}/
(\lambda t^{\nu /2})\cos ( \nu \pi /2) }\nonumber\\[-8pt]\\[-8pt]
&&\phantom{\frac{1}{\pi \nu |x|}\int_{0}^{+\infty }}
{}\times\sin \biggl( \frac{%
|x|w^{\nu /2}}{\lambda t^{\nu /2}}\sin \biggl(\frac{\nu \pi }{2}\biggr)\biggr) \,dw,
\nonumber
\end{eqnarray}%
for $0<\nu <2$.
\end{corollary}

\begin{pf}
The first integration in (\ref{sin}) gives
\begin{eqnarray*}
u_{\nu }(x,t)
&=&\frac{1}{\pi \nu |x|\sin (\nu \pi /2)} \cos \biggl( \frac{%
\nu \pi }{2}-\frac{|x|w^{\nu /2}}{\lambda t^{\nu /2}}\sin \biggl(\frac{\nu \pi }{2}%
\biggr)\biggr)\\
&&{}\times e^{-w}e^{-|x|w^{\nu /2}/(\lambda t^{\nu /2})\cos (
\nu \pi /2) }\Bigm\vert _{0}^{+\infty } \\
&&{}+\frac{1}{\pi \nu |x|\sin (\nu \pi /2)}\int_{0}^{+\infty
}e^{-w}e^{-|x|w^{\nu /2}/(\lambda t^{\nu /2})\cos ( \nu \pi
/2) }\\
&&\phantom{{}+\frac{1}{\pi \nu |x|\sin (\nu \pi /2)}\int_{0}^{+\infty}}
{}\times\cos \biggl( \frac{\nu \pi }{2}-\frac{|x|w^{\nu /2}}{\lambda
t^{\nu /2}}\sin \biggl(\frac{\nu \pi }{2}\biggr)\biggr) \,dw \\
&&{}+\frac{\cos ( \nu \pi /2) }{2\pi \sin (\nu \pi /2%
)\lambda t^{\nu /2}}\int_{0}^{+\infty }e^{-w}w^{\nu /2-1}e^{-
|x|w^{\nu /2}/(\lambda t^{\nu /2})\cos ( \nu \pi /2)
}\\
&&\phantom{{}+\frac{\cos ( \nu \pi /2) }{2\pi \sin (\nu \pi /2%
)\lambda t^{\nu /2}}\int_{0}^{+\infty }}
{}\times\cos \biggl( \frac{\nu \pi }{2}-\frac{|x|w^{\nu /2}}{\lambda t^{\nu /2}}\sin
\biggl(\frac{\nu \pi }{2}\biggr)\biggr)\, dw
\\
&=&-\frac{\cot (\nu \pi /2)}{\pi \nu |x|}\\
&&{}+\frac{1}{\pi \nu |x|\sin (%
\nu \pi /2)}\int_{0}^{+\infty }e^{-w}e^{-|x|w^{\nu /2}/
(\lambda t^{\nu /2})\cos ( \nu \pi /2) }\\
&&\phantom{{}+\frac{1}{\pi \nu |x|\sin (%
\nu \pi /2)}\int_{0}^{+\infty }}
{}\times\cos \biggl( \frac{%
\nu \pi }{2}-\frac{|x|w^{\nu /2}}{\lambda t^{\nu /2}}\sin \biggl(\frac{\nu \pi }{2}%
\biggr)\biggr)\, dw \\
&&{}- \frac{\cos ( \nu \pi /2) }{\pi \nu |x|\sin
^{2}( \nu \pi /2) }e^{-w}e^{-|x|w^{\nu /2}/(\lambda
t^{\nu /2})\cos ( \nu \pi /2) }\\
&&\quad{}\times\sin \biggl( \frac{\nu \pi
}{2}-\frac{|x|w^{\nu /2}}{\lambda t^{\nu /2}}\sin \biggl(\frac{\nu \pi }{2}%
\biggr)\biggr) \Bigm\vert _{0}^{+\infty } \\
&&{}-\frac{\cos ( \nu \pi /2) }{\pi \nu |x|\sin ^{2}(
\nu \pi /2)}\int_{0}^{+\infty }e^{-w}e^{-|x|w^{\nu /2}/(\lambda
t^{\nu /2})\cos ( \nu \pi /2) }\\
&&\phantom{{}-\frac{\cos ( \nu \pi /2) }{\pi \nu |x|\sin ^{2}(
\nu \pi /2)}\int_{0}^{+\infty }}
{}\times\sin \biggl( \frac{\nu \pi
}{2}-\frac{|x|w^{\nu /2}}{\lambda t^{\nu /2}}\sin \biggl(\frac{\nu \pi }{2}%
\biggr)\biggr) \,dw \\
&&{}-\frac{\cos ^{2}( \nu \pi /2) }{2\pi \sin ^{2}(
\nu \pi /2)\lambda t^{\nu /2}}\int_{0}^{+\infty }e^{-w}w^{\nu /2
-1}e^{-|x|w^{\nu /2}/(\lambda t^{\nu /2})\cos ( \nu \pi /2%
) }\\
&&\phantom{-\frac{\cos ^{2}( \nu \pi /2) }{2\pi \sin ^{2}(
\nu \pi /2)\lambda t^{\nu /2}}\int_{0}^{+\infty }}
{}\times\sin \biggl( \frac{\nu \pi }{2}-\frac{|x|w^{\nu /2}}{\lambda t^{\nu
/2}}\sin \biggl(\frac{\nu \pi }{2}\biggr)\biggr)\, dw.
\end{eqnarray*}

Therefore, from (\ref{sin}) we have that%
\begin{eqnarray*}
&&\biggl\{ 1+\frac{\cos ^{2}( \nu \pi /2) }{\sin ^{2}(%
\nu \pi /2)}\biggr\} u_{\nu }(x,t) \\
&&\qquad=\frac{1}{\pi \nu |x|\sin (\nu \pi /2)}\\
&&\quad\qquad
{}\times\int_{0}^{+\infty
}e^{-w}e^{-|x|w^{\nu /2}/(\lambda t^{\nu /2})\cos ( \nu \pi
/2) }\\
&&\quad\qquad\phantom{{}\times\int_{0}^{+\infty}}{}\times \biggl[ \cos \biggl( \frac{\nu \pi }{2}-\frac{|x|w^{\nu /2}}{\lambda
t^{\nu /2}}\sin \biggl(\frac{\nu \pi }{2}\biggr)\biggr)\\
&&\hspace*{90pt}
{} -\cot \biggl(\frac{\nu \pi }{2}\biggr)\sin
\biggl( \frac{\nu \pi }{2}-\frac{|x|w^{\nu /2}}{\lambda t^{\nu /2}}\sin \biggl(%
\frac{\nu \pi }{2}\biggr)\biggr) \biggr]\, dw
\\
&&\qquad=\frac{1}{\pi \nu |x|\sin (\nu \pi /2)}\int_{0}^{+\infty
}e^{-w}e^{-|x|w^{\nu /2}/(\lambda t^{\nu /2})\cos ( \nu \pi
/2) } \\
&&\quad\qquad{}\times \biggl\{\biggl[ \cos \biggl(\frac{\nu \pi }{2}\biggr)\cos \biggl( \frac{|x|w^{\nu /2}}{%
\lambda t^{\nu /2}}\sin \biggl(\frac{\nu \pi }{2}\biggr)\biggr)\\
&&\hspace*{56pt}
{} +\sin \biggl(\frac{\nu \pi }{2}%
\biggr)\sin \biggl( \frac{|x|w^{\nu /2}}{\lambda t^{\nu /2}}\sin \biggl(\frac{\nu \pi }{2}%
\biggr)\biggr) \biggr] \\
&&\quad\qquad\hspace*{17pt}{}-\cot \biggl(\frac{\nu \pi }{2}\biggr)\biggl[ \sin \biggl(\frac{\nu \pi }{2}\biggr)\cos \biggl(
\frac{|x|w^{\nu /2}}{\lambda t^{\nu /2}}\sin \biggl(\frac{\nu \pi
}{2}\biggr)\biggr)\\
&&\hspace*{110pt}
{}-\cos \biggl(\frac{\nu \pi }{2}\biggr)\sin \biggl( \frac{|x|w^{\nu /2}}{\lambda t^{\nu /2}%
}\sin \biggl(\frac{\nu \pi }{2}\biggr)\biggr) \biggr] ,
\end{eqnarray*}%
which easily gives (\ref{sin.3}).
\end{pf}

\begin{remark}\label{rem5.3}
We can check that, for $\nu =1$, (\ref{sin.3}) reduces again to the Gaussian
density:
\begin{eqnarray*}
u_{1}(x,t)
&=&\frac{1}{\pi |x|}\int_{0}^{+\infty }e^{-w}\sin \biggl( \frac{|x|w^{1/2}}{%
\lambda t^{1/2}}\biggr) \, dw \\
&=&\biggl[ w=\frac{y^{2}}{2}\frac{\lambda ^{2}t}{|x|^{2}}\biggr] \\
&=&\frac{\lambda ^{2}t}{\pi |x|^{3}}\int_{0}^{+\infty }ye^{-y^{2}/2
\lambda ^{2}t/|x|^{2}}\sin \frac{y}{\sqrt{2}}\,dy \\
&=&\frac{1}{\sqrt{2}\pi |x|}\int_{0}^{+\infty }e^{-y^{2}/2
\lambda ^{2}t/|x|^{2}}\cos \frac{y}{\sqrt{2}}\,dy \\
&=&\frac{1}{2\sqrt{\pi t\lambda ^{2}}}e^{-x^{2}/(4t\lambda ^{2})}
\end{eqnarray*}%
as in (\ref{gaus}).

With respect to (\ref{sin}), formula (\ref{sin.3}) is more appealing as it
allows an easier analysis of the limit for $|x|\rightarrow 0$:%
\begin{eqnarray}
\lim_{|x|\rightarrow 0}u_{\nu }(x,t) &=&\frac{1}{\pi \nu }\frac{\sin (
\nu \pi /2)}{\lambda t^{\nu /2}}\int_{0}^{+\infty }w^{\nu /2}e^{-w}\,dw
\label{lim} \\
&=&\frac{1}{\pi \nu }\frac{\sin (\nu \pi /2)}{\lambda t^{\nu /2}}%
\Gamma \biggl( \frac{\nu }{2}+1\biggr) .  \nonumber
\end{eqnarray}

For $t\rightarrow +\infty $, (\ref{lim}) decreases for all values of $\nu \in
( 0,2] .$

Moreover in the case $\nu =1$, formula (\ref{lim}) gives the maximum value of
the Brownian density. For $\nu =2$ (\ref{lim}) is zero for all $t>0$,
because in this case (\ref{equation1}) becomes the wave equation and its
solution has the form of the sum of Dirac's impulse functions travelling in
opposite directions.

By means of the following formula%
\[
\int_{0}^{+\infty }\frac{\sin qx}{x}e^{-px}\,dx=\arctan \frac{q}{p},\qquad p>0
\]
[Gradshteyn and Ryzhik (\citeyear{GraRyz1994}), formula 3.941.1, page 523] we can check that (%
\ref{sin.3}) integrates to one, as follows:%
\begin{eqnarray*}
\int_{-\infty }^{+\infty }u_{\nu }(x,t)\,dx
&=&\frac{2}{\pi \nu }\int_{0}^{+\infty }e^{-w}\int_{0}^{+\infty }\frac{1}{x}%
e^{-xw^{\nu /2}/(\lambda t^{\nu /2})\cos ( \nu \pi /2
) }\\
&&\hspace*{88pt}
{}\times\sin \biggl( \frac{xw^{\nu /2}}{\lambda t^{\nu /2}}\sin \biggl(\frac{\nu
\pi }{2}\biggr)\biggr)\, dx\,dw \\
&=&\frac{2}{\pi \nu }\frac{\nu \pi }{2}\int_{0}^{+\infty }e^{-w}\,dw=1.
\end{eqnarray*}

\begin{figure}

\includegraphics{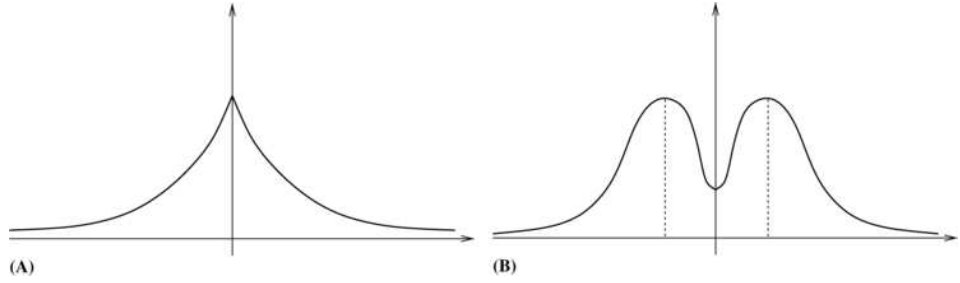}

\caption{(\textup{A}) The solution $u_{\nu}(x,t)$, for $0<\nu<1$.  (\textup{B}) The
solution $u_{\nu}(x,t)$, for $1<\nu<2$.}\label{fig2}
\end{figure}

Finally it is interesting to analyze the behavior of the solution (for $x$
varying and $t$ fixed), which is substantially different in the two
intervals $0<\nu \leq 1$ and $1<\nu \leq 2$ (see Figure \ref{fig2} above). We
rewrite formula (\ref{sin.3}) as follows: for $x>0,$%
\[
u_{\nu }(x,t)=\frac{1}{\pi \nu }\int_{0}^{+\infty }\frac{g(x,w,t)}{x}%
e^{-w}\,dw,
\]
where $g(x,w,t)=e^{-xA\cos (\nu \pi /2)}\sin ( xA\sin (\nu
\pi /2)) $ and $A=w^{\nu /2}/\lambda t^{\nu /2}$.

The first derivative of $\frac{g(x,w,t)}{x}$ with respect to $x$ is equal to
zero if%
\begin{equation}
\frac{g_{x}}{g}=\frac{1}{x},  \label{eq}
\end{equation}%
where%
\begin{eqnarray*}
g_{x} &=&-A\cos \frac{\nu \pi }{2}e^{-xA\cos (\nu \pi /2)}\sin \biggl(
xA\sin \frac{\nu \pi }{2}\biggr)  \\
&&{}+A\sin \frac{\nu \pi }{2}e^{-xA\cos (\nu \pi /2)}\cos \biggl( xA\sin
\frac{\nu \pi }{2}\biggr) \\
&=&Ae^{-xA\cos (\nu \pi /2)}\sin \biggl( \frac{\nu \pi }{2}-xA\sin
\frac{\nu \pi }{2}\biggr) .
\end{eqnarray*}

The solution to (\ref{eq}) is%
\[
\lg g=\lg x+const
\]
or, otherwise,
\[
g=xconst.
\]
By choosing $const=1$, we obtain that $u_{\nu }(x,t)$ attains its maximum on
the positive half-line if
\begin{equation}
xe^{xA\cos (\nu \pi /2)}=\sin \biggl( xA\sin \frac{\nu \pi }{2}\biggr) .
\label{max}
\end{equation}

For $1<\nu \leq 2$ there exists only one value of $x$ which verifies the
condition (\ref{max}) and this is in accordance with the behavior of the
solutions $u_{\nu }$ presented in Fujita (\citeyear{Fuj1990}), where the relationship with
stable laws is exploited.

On the other hand, for $0<\nu \leq 1,$ no positive value satisfies (\ref{max}%
) and therefore the maximum is in the origin. The previous results are
confirmed by the following theorems.
\end{remark}

We now present the general results concerning the relationship between the
solution $u_{\nu }(x,t)$ and the stable densities. We need to analyze the
two intervals $0<\nu \leq 1$ and $1<\nu \leq 2$ separately.

\begin{thm}\label{thm5.4}
For $0<\nu \leq 1$, the solution
to
\begin{equation}
\cases{
\displaystyle\frac{\partial ^{\nu }u}{\partial t^{\nu }}=\lambda ^{2}\frac{\partial ^{2}u%
}{\partial x^{2}}, \cr
u(x,0)=\delta (x),
}
\qquad x\in \mathbb{R}, t>0,  \label{stable}
\end{equation}%
can be represented as
\begin{eqnarray} \label{stable.1}
u_{\nu }(x,t) &=&\frac{1}{\nu }\frac{\lambda ^{2/\nu }t}{|x|^{2/\nu
+1}}p_{\nu /2}\biggl( \frac{\lambda ^{2/\nu }t}{|x|^{2/\nu }};
\frac{\nu }{2},1\biggr) \nonumber\\[-8pt]\\[-8pt]
&=&\frac{1}{\nu |x|^{2/\nu +1}}p_{\nu /2}\biggl( \frac{1}{%
|x|^{2/\nu }};\frac{\nu }{2},\frac{1}{\lambda t^{\nu /2}}%
\biggr)  \nonumber
\end{eqnarray}%
where $p_{\frac{\nu }{2}}( \cdot ;\frac{\nu }{2},1) $
is the density of a stable distribution of parameters $\gamma =\frac{\nu
}{2}$ and $\eta =1$; see (\ref{stab}).
\end{thm}

\begin{pf}
From (\ref{equation2}), by using the reflection
formula for the Gamma function we have that%
\begin{eqnarray} \label{stable.2}
u_{\nu }(x,t) &=&\frac{1}{2\lambda t^{\nu /2}}\sum_{k=0}^{\infty }\biggl( -%
\frac{|x|}{\lambda t^{\nu /2}}\biggr) ^{k}\frac{1}{k!\Gamma (1-\nu
(k+1)/2)}\nonumber\\[-8pt]\\[-8pt]
&=&\frac{1}{2\lambda t^{\nu /2}}\sum_{k=0}^{\infty }\biggl( -\frac{|x|}{%
\lambda t^{\nu /2}}\biggr) ^{k}\frac{1}{k!}\frac{\sin (\pi \nu (k+1)/2)%
}{\pi }\Gamma \biggl( \frac{\nu (k+1)}{2}\biggr) .  \nonumber
\end{eqnarray}

In view of the series representation of stable functions, which for $%
0<\alpha <1$ reads%
\[
p_{\alpha }(x;\gamma ,1)=\frac{\alpha }{\pi }\sum_{r=0}^{\infty }(-1)^{r}%
\frac{\Gamma (\alpha (r+1))}{r!}x^{-\alpha (r+1)-1}\sin \biggl[ \frac{\pi }{2}%
(\gamma +\alpha )(r+1)\biggr]
\]
[see Feller (\citeyear{Fel1971}), formula (6.10), page 583, with some corrections,  Lukacs
(\citeyear{Luk1969}) and Zolotarev (\citeyear{Zol1986})], we can obtain the first expression in (\ref{stable.1}). The second
expression can be derived by applying the self-similarity property of the
stable random variables.
\end{pf}

Finally we consider the case $1\leq \nu \leq 2$ and we state the
following result:

\begin{thm}\label{thm5.5}
The solution to
\begin{equation}
\cases{
\displaystyle\frac{\partial ^{\nu }u}{\partial t^{\nu }}=\lambda ^{2}\frac{\partial ^{2}u%
}{\partial x^{2}}, \cr
u(x,0)=\delta (x), \cr
u_{t}(x,0)=0,}%
\qquad x\in \mathbb{R}, t>0,  \label{gen}
\end{equation}%
for $1\leq \nu \leq 2$, can be represented as
\begin{eqnarray}\label{gen.2}
u_{\nu }(x,t) &=&\frac{2}{\nu }\frac{1}{2\lambda t^{\nu /2}}p_{2/\nu
}\biggl( \frac{|x|}{\lambda t^{\nu /2}};\frac{2}{\nu }(\nu -1),1\biggr)
\nonumber\\[-8pt]\\[-8pt]
&=&\frac{1}{\nu }p_{2/\nu }\biggl( |x|;\frac{2}{\nu }(\nu -1),\lambda
^{2/\nu }t\biggr) ,  \nonumber
\end{eqnarray}
where $p_{2/\nu }(\cdot ;\frac{2}{\nu }(\nu -1),1)$ is
the density of a stable distribution of parameters $\gamma =\frac{2}{\nu }%
(\nu -1)$ and $\eta =1$.
\end{thm}

\begin{pf}
By following the same steps as in the previous
theorem we can recognize in (\ref{stable.2}), up to the normalizing
constant, the series representation of the stable law $p_{2/\nu }$
of order $\alpha =2/\nu $ [see (\ref{qua.7})], so that we get (\ref{gen.2}).
\end{pf}

\begin{remark}\label{rem5.4}
In view of Theorems \ref{thm2.3} and \ref{thm5.5} and by considering the property of
self-similarity of the stable laws, we can write that%
\begin{eqnarray}\label{stable.7}
u_{\nu }(x,t) &=&\frac{1}{\nu }\int_{0}^{\infty }\frac{e^{-x^{2}/
(4w\lambda)}}{\sqrt{4\pi w\lambda }}\frac{1}{\lambda t^{\nu }}p_{1/
\nu }\biggl( \frac{|w|}{\lambda t^{\nu }},\frac{1}{\nu }(2\nu -1),1\biggr)\, dw
\nonumber\\[-8pt]\\[-8pt]
&=&\frac{1}{\nu }\int_{0}^{\infty }\frac{e^{-x^{2}/(4w\lambda)}}{%
\sqrt{4\pi w\lambda }}p_{1/\nu }\biggl( |w|,\frac{1}{\nu }(2\nu
-1),\lambda ^{1/\nu }t\biggr)\, dw.  \nonumber
\end{eqnarray}

Formula (\ref{stable.7}) shows that the solution $u_{\nu },$ for $\frac{1}{2}%
<\nu \leq 1,$ can be interpreted as the distribution of the process $B(|%
\mathcal{S}_{\nu }(t)|),$ $t>0$, where $\mathcal{S}_{\nu }$ is the stable
process with density $\frac{1}{\nu }p_{1/\nu }( |\cdot |,\frac{1%
}{\nu }(2\nu -1),\lambda ^{1/\nu }t) $.

Moreover, as a consequence of Theorems 2.1 and 5.5, the solution of our
problem (\ref{equation1})--(\ref{condition 1}), for $\frac{1}{2}<\nu \leq 1$,
can be written in an alternative to the form (\ref{vologda}) also as a stable
law evaluated at a Brownian time:%
\[
u_{\nu }(x,t)=\frac{1}{\nu }\int_{0}^{\infty }\frac{e^{-s^{2}/(4t)}}{%
\sqrt{\pi t}}\frac{1}{2\lambda s^{\nu }}p_{1/\nu }\biggl( \frac{|x|}{%
\lambda s^{\nu }},\frac{1}{\nu }(2\nu -1),1\biggr) \,ds.
\]
\end{remark}

\begin{remark}\label{rem5.5}
We check that, for $\nu =1$, both expressions (\ref{stable.1}) and (\ref%
{gen.2}) yield the Gaussian density%
\begin{equation}
u_{1}(x,t)=\frac{1}{2\lambda \sqrt{\pi t}}e^{-x^{2}/(4\lambda ^{2}t)}.
\label{gen.6}
\end{equation}

We start by considering the last expression in (\ref{stable.1})%
\begin{equation}
u_{1}(x,t)=\frac{1}{|x|^{3}}p_{1/2}\biggl( \frac{1}{|x|^{2}};\frac{1}{%
2},\frac{1}{\lambda t^{1/2}}\biggr) ,  \label{gen.7}
\end{equation}%
where [from (\ref{stab})], for $y>0,$%
\begin{eqnarray} \label{gen.8}
p_{1/2}\biggl( y;\frac{1}{2},\frac{1}{\lambda t^{1/2}}\biggr)
&=&\frac{1}{2\pi }\int_{-\infty }^{+\infty }e^{-i\beta y}\exp \biggl\{ -\frac{%
|\beta |^{1/2}}{\lambda t^{1/2}}e^{-i\pi /4 \beta /|\beta |
}\biggr\} \,d\beta \nonumber\\[-8pt]\\[-8pt]
&=&\frac{1}{\sqrt{2}\lambda t^{1/2}}\frac{e^{-1/(2y(\sqrt{2}\lambda
t^{1/2})^{2})}}{\sqrt{2\pi y^{3}}}.  \nonumber
\end{eqnarray}

By taking in (\ref{gen.8}) $y=\frac{1}{|x|^{2}}$ we get from (\ref{gen.7}%
) the Gaussian density (\ref{gen.6}).
Formula (\ref{gen.2}) immediately supplies (\ref{gen.6}) for $\nu =1.$
\end{remark}

\section*{Acknowledgments}

The authors thank one anonymous referee for bringing their attention to some
relevant papers on fractional equations. Thanks are also due for his
accurate check of the text and of the calculations.

\def\cprime{$'$}

\printaddresses

\end{document}